\documentclass{amsart} 
\usepackage{latexsym} 
\usepackage{amssymb} 
\usepackage{amsmath}
\usepackage{fancyhdr}
\begin{document}


\newtheorem{theorem}{Theorem} 
\newtheorem{problem}{Problem} 
\newtheorem{definition}{Definition} 
\newtheorem{lemma}{Lemma} 
\newtheorem{proposition}{Proposition} 
\newtheorem{corollary}{Corollary} 
\newtheorem{example}{Example} 
\newtheorem{conjecture}{Conjecture} 
\newtheorem{algorithm}{Algorithm} 
\newtheorem{exercise}{Exercise} 
\newtheorem{remarkk}{Remark} 
 
\newcommand{\be}{\begin{equation}} 
\newcommand{\ee}{\end{equation}} 
\newcommand{\bea}{\begin{eqnarray}} 
\newcommand{\eea}{\end{eqnarray}} 
\newcommand{\beq}[1]{\begin{equation}\label{#1}} 
\newcommand{\eeq}{\end{equation}} 
\newcommand{\beqn}[1]{\begin{eqnarray}\label{#1}} 
\newcommand{\eeqn}{\end{eqnarray}} 
\newcommand{\beaa}{\begin{eqnarray*}} 
\newcommand{\eeaa}{\end{eqnarray*}} 
\newcommand{\req}[1]{(\ref{#1})} 
 
\newcommand{\lip}{\langle} 
\newcommand{\rip}{\rangle} 

\newcommand{\uu}{\underline} 
\newcommand{\oo}{\overline} 
\newcommand{\La}{\Lambda} 
\newcommand{\la}{\lambda} 
\newcommand{\eps}{\varepsilon} 
\newcommand{\om}{\omega} 
\newcommand{\Om}{\Omega} 
\newcommand{\ga}{\gamma} 
\newcommand{\rrr}{{\Bigr)}} 
\newcommand{\qqq}{{\Bigl\|}} 
 
\newcommand{\dint}{\displaystyle\int} 
\newcommand{\dsum}{\displaystyle\sum} 
\newcommand{\dfr}{\displaystyle\frac} 
\newcommand{\bige}{\mbox{\Large\it e}} 
\newcommand{\integers}{{\Bbb Z}} 
\newcommand{\rationals}{{\Bbb Q}} 
\newcommand{\reals}{{\rm I\!R}} 
\newcommand{\realsd}{\reals^d} 
\newcommand{\realsn}{\reals^n} 
\newcommand{\NN}{{\rm I\!N}} 
\newcommand{\DD}{{\rm I\!D}} 
\newcommand{\M}{{\rm I\!M}} 
\newcommand{\degree}{{\scriptscriptstyle \circ }} 
\newcommand{\dfn}{\stackrel{\triangle}{=}} 
\def\complex{\mathop{\raise .45ex\hbox{${\bf\scriptstyle{|}}$} 
     \kern -0.40em {\rm \textstyle{C}}}\nolimits} 
\def\hilbert{\mathop{\raise .21ex\hbox{$\bigcirc$}}\kern -1.005em {\rm\textstyle{H}}} 
\newcommand{\RAISE}{{\:\raisebox{.6ex}{$\scriptstyle{>}$}\raisebox{-.3ex} 
           {$\scriptstyle{\!\!\!\!\!<}\:$}}} 
 
\newcommand{\hh}{{\:\raisebox{1.8ex}{$\scriptstyle{\degree}$}\raisebox{.0ex} 
           {$\textstyle{\!\!\!\! H}$}}} 

\newcommand{\OO}{\won} 
\newcommand{\calA}{{\mathcal A}} 
\newcommand{\calB}{{\mathcal B}} 
\newcommand{\calC}{{\cal C}} 
\newcommand{\calD}{{\cal D}} 
\newcommand{\calE}{{\cal E}} 
\newcommand{\calF}{{\mathcal F}} 
\newcommand{\calG}{{\cal G}} 
\newcommand{\calH}{{\cal H}} 
\newcommand{\calK}{{\cal K}} 
\newcommand{\calL}{{\mathcal L}} 
\newcommand{\calM}{{\cal M}} 
\newcommand{\calO}{{\cal O}} 
\newcommand{\calP}{{\cal P}} 
\newcommand{\calS}{{\mathcal S}} 
\newcommand{\calU}{{\mathcal U}} 
\newcommand{\calX}{{\cal X}} 
\newcommand{\calXX}{{\cal X\mbox{\raisebox{.3ex}{$\!\!\!\!\!-$}}}} 
\newcommand{\calXXX}{{\cal X\!\!\!\!\!-}} 
\newcommand{\gi}{{\raisebox{.0ex}{$\scriptscriptstyle{\cal X}$} 
\raisebox{.1ex} {$\scriptstyle{\!\!\!\!-}\:$}}} 
\newcommand{\hna}{\hat{\nabla}}
\newcommand{\intsim}{\int_0^1\!\!\!\!\!\!\!\!\!\sim} 
\newcommand{\intsimt}{\int_0^t\!\!\!\!\!\!\!\!\!\sim} 
\newcommand{\pp}{{\partial}} 
\newcommand{\al}{{\alpha}} 
\newcommand{\sB}{{\cal B}} 
\newcommand{\sL}{{\cal L}} 
\newcommand{\sF}{{\cal F}} 
\newcommand{\sE}{{\cal E}} 
\newcommand{\sX}{{\cal X}} 
\newcommand{\R}{{\rm I\!R}} 
\renewcommand{\L}{{\rm I\!L}} 
\newcommand{\vp}{\varphi} 
\newcommand{\N}{{\rm I\!N}} 
\def\ooo{\lip} 
\def\ccc{\rip} 
\newcommand{\ot}{\hat\otimes} 
\newcommand{\rP}{{\Bbb P}} 
\newcommand{\bfcdot}{{\mbox{\boldmath$\cdot$}}} 
 
\renewcommand{\varrho}{{\ell}} 
\newcommand{\dett}{{\textstyle{\det_2}}} 
\newcommand{\sign}{{\mbox{\rm sign}}} 
\newcommand{\TE}{{\rm TE}} 
\newcommand{\TA}{{\rm TA}} 
\newcommand{\E}{{\rm E\,}} 
\newcommand{\won}{{\mbox{\bf 1}}} 
\newcommand{\Lebn}{{\rm Leb}_n} 
\newcommand{\Prob}{{\rm Prob\,}} 
\newcommand{\sinc}{{\rm sinc\,}} 
\newcommand{\ctg}{{\rm ctg\,}} 
\newcommand{\loc}{{\rm loc}} 
\newcommand{\trace}{{\,\,\rm trace\,\,}} 
\newcommand{\Dom}{{\rm Dom}} 
\newcommand{\ifff}{\mbox{\ if and only if\ }} 
\newcommand{\nproof}{\noindent {\bf Proof:\ }} 
\newcommand{\remark}{\noindent {\bf Remark:\ }} 
\newcommand{\remarks}{\noindent {\bf Remarks:\ }} 
\newcommand{\note}{\noindent {\bf Note:\ }}

\newcommand{\boldx}{{\bf x}} 
\newcommand{\boldX}{{\bf X}} 
\newcommand{\boldy}{{\bf y}} 
\newcommand{\boldR}{{\bf R}} 
\newcommand{\uux}{\uu{x}} 
\newcommand{\uuY}{\uu{Y}} 
 
\newcommand{\limn}{\lim_{n \rightarrow \infty}} 
\newcommand{\limN}{\lim_{N \rightarrow \infty}} 
\newcommand{\limr}{\lim_{r \rightarrow \infty}} 
\newcommand{\limd}{\lim_{\delta \rightarrow \infty}} 
\newcommand{\limM}{\lim_{M \rightarrow \infty}} 
\newcommand{\limsupn}{\limsup_{n \rightarrow \infty}} 
 
\newcommand{\ra}{ \rightarrow }

\newcommand{\ARROW}[1] 
  {\begin{array}[t]{c}  \longrightarrow \\[-0.2cm] \textstyle{#1} \end{array} } 
 
\newcommand{\AR} 
 {\begin{array}[t]{c} 
  \longrightarrow \\[-0.3cm] 
  \scriptstyle {n\rightarrow \infty} 
  \end{array}} 
 
\newcommand{\pile}[2] 
  {\left( \begin{array}{c}  {#1}\\[-0.2cm] {#2} \end{array} \right) } 
 
\newcommand{\floor}[1]{\left\lfloor #1 \right\rfloor} 
 
\newcommand{\mmbox}[1]{\mbox{\scriptsize{#1}}} 
 
\newcommand{\ffrac}[2] 
  {\left( \frac{#1}{#2} \right)} 
 
\newcommand{\one}{\frac{1}{n}\:} 
\newcommand{\half}{\frac{1}{2}\:} 

\def\le{\leq} 
\def\ge{\geq} 
\def\lt{<} 
\def\gt{>} 
 
\def\squarebox#1{\hbox to #1{\hfill\vbox to #1{\vfill}}} 
\newcommand{\nqed}{\hspace*{\fill} 
           \vbox{\hrule\hbox{\vrule\squarebox{.667em}\vrule}\hrule}\bigskip}  
 
\title{Variational Solutions of the Monge Transport Problem and the
  Monge-Amp\`ere Equation in Abstract Wiener Space}

\author{ A. S. \"Ust\"unel} 
\maketitle

\noindent 
{\bf Abstract:}{\small{ Let $(W,H,\mu)$ be an abstract  Wiener space,
    assume that $T=I_W+\nabla\varphi$ is the solution of the Monge
    problem associated to the measures $d\mu$ and
    $d\nu=Ld\mu=e{^-f}d\mu$}. Under the finite information hypothesis, using
  a variational method, we
  prove that the forward potential of the Monge-Kantorovitch problem
  satisfies the ``structure equation''
  $$
\delta((I_H+\nabla^2\varphi)^{-1}-I_H)=\nabla\varphi+\nabla f\circ
  T
$$
 and with it the Sobolev regularity of the backward
  Monge potential is proven map. A similar structure equations  also holds for the forward
  Monge potential and it implies the regularity of it for $1-\eps$
  log-concave densities. We show that $L=e^{-f}$ can be represented as 
$$
L=\dett(I+\nabla^2\psi)\exp\left[-\calL\psi-\half|\nabla\psi|_H^2\right]\,,
$$
where $\psi$ is the backward Monge potential. Moreover the forward
potentail gives the solution of the Monge-Amp\`ere equation:
$$
L\circ T\,\La=1
$$
$\mu$-a.s., where 
$$
\La=\dett(I+\nabla^2_a\varphi)\left[-\calL^a\varphi-\half|\varphi|_H^2\right]\,,
$$
and $\nabla^2_a\varphi$ is the Radon-Nikodym derivative of the
Hilbert-Schmidt-valued absolutely continuous part of the vector
measure $\nabla^2\varphi$. In particular, for $L\neq 0$ almost surely,
the Girsanov identity holds:
$$
\int_W g\circ T \La(\varphi)d\mu=\int_W gd\mu\,,
$$
for any $g\in C_b(W)$.
}
\vspace{0.5cm} 
\tableofcontents

\section{\bf{Introductions} }
\noindent
Let $\nu$ be the probability measure defined by
\begin{equation}
\label{target}
d\nu=\frac{1}{c} e^{-f}d\mu
\end{equation}
such that the relative entropy of $\nu$ w.r.t. the Wiener measure $\mu$, denoted as
$H(\nu|\mu)$ is finite. Let $\Sigma(\mu,\nu)$ be the set of the
probability measures on $(W\times W,\calB(W\times W))$ whose first
marginals are $\mu$ and the secones ones are $\nu$. Consider the
problem of minimization which defines also a strong Wasserstein distance
between $\mu$ and $\nu$:
$$
\inf\left(\int_{W\times W}|x-y|_H^2d\beta (x,y):\,\beta\in
  \Sigma(\mu,\nu)\right)=d_2^2(\mu,\nu)\,,
$$
where $|\cdot|_H$ denotes the Cameron-Martin norm. In the finite
dimensional case this problem has been extensively studied since
almost  three centuries and we refer  to the texts
\cite{R-R} and \cite{Vil} for history and references and also to
\cite{BRE} and \cite{Mc}.

In the infinite dimensional
case, where the cost function is very singular, in the sense that the
set on which the cost function is finite has zero measure w.r.t. the
product measure $\mu\times \nu$ has
been solved in a series of papers (\cite{fandu1, fandu2,fandu3}) and
the answer can be summarized as follows: There exists a $1$-convex
function $\varphi:W\to\R$, in the Gaussian Sobolev space $\DD_{2,1}$,
called Monge potential or Monge-Brenier map such that the above
infimum is attained at $\ga=(I_W\times T)\mu$, i.e., the image of the
measure $\mu$ under the map $I_W\times T$, where $T=I_W+\nabla
\varphi$, where $\nabla\varphi$ is the $L^2(\mu)$-extended derivative
of $\varphi$ in the direction of Cameron-Martin space. Moreover, there
exists also a dual Monge potential $\psi:W\to \R$, which has an
$L^2(\nu)$-extended derivative in the direction of Cameron-Martin
space, such that, the map $S=I_W+\nabla\psi$ satisfies $(S\times
I_W)\nu=(I_W\times T)\mu=\ga$, hence $T\circ S=I_W$ $\nu$-a.s. and
$S\circ T=I_W$ $\mu$-a.s. The next important issue in this subject is
to show the Sobolev regularity of the Monge-Brenier maps in such a way
that one can write the Jacobian functions associated to the
corresponding transformations $T$ and $S$. In finite dimensional case
this problem has been treated by several authors (cf. \cite{Caf} and
the references given in \cite{Vil}). In the infinite dimensional case
there are also some results (cf. \cite{Kol, F-N, fandu4} ) which are
generalizations of the results given in \cite{fandu2,fandu3}. These
results are generally some extensions of the methods developped
especially by L. Caffarelli, though we have also  given
another method to calculate the Jacobian functions in infinite
dimensions using the It\^o calculus.

 In this work we shall present a
totally different method, namely, we shall prove the Sobolev
regularity of the Monge-Brenier functions using the Calculus of
variations. Let us begin by recalling a celebrated variational
formula, which holds on any measurable space but we formulate it on a
Wiener space for the notational simplicity:
\begin{equation}
\label{1st-variation}
-\log\int_W e^{-f}d\mu=\inf\left(\int_W fd\ga+H(\ga|\mu):\,\ga\in M_1(W)\right)
\end{equation}
where $M_1(W)$ denotes the set of probability measures on $(W,\calF)$,
$\calF$ being the Borel sigma field of $W$, $\ga,\,\mu$ are as
described above and suppose that the measure $e^{-f}d\mu$  is of
finite relative entropy w.r.t. $\mu$.  Then the
infimum is attained at $d\nu=e^{-f}d\mu$ provided that $H(\nu|\mu)$ is finite,
cf. \cite{ASU-3}, \cite{ASU-4}. On the other hand, we know from \cite{fandu1} that
there exists some $\varphi\in\DD_{2,1}$, $1$-convex function such that
$(I_W+\nabla\varphi)\mu=\nu$, where we use the same notation for the
image of a point and of a measure under a measurable map (here the map
under question is  $T=I_W+\nabla\varphi$). Consequently the following
identity holds true:
$$
-\log\int e^{-f}d\mu=\inf\left(\int f\circ M
  d\mu+H(M\mu|\mu):\,M=I_W+\nabla a,\,a\in \DD_{2,1}\right)\,.
$$
Therefore 
\begin{equation}
\label{ineq_1}
-\log\int e^{-f}d\mu\geq \inf\left(\int f\circ (I_W+\xi)
  d\mu+H((I_W+\xi)\mu|\mu):\,\xi\in \DD_{2,0}(H)\right)\,.
\end{equation}
For this infimum to be finite we need that
$H((I_W+\xi)\mu|\mu)<\infty$, which implies 
$(I_W+\xi)\mu\ll\mu$. Besides the right hand side of the inequality
(\ref{ineq_1}) is always greater than
$$
\inf\left(\int fd\ga+H(\ga|\mu):\,\ga\in M_1(W)\right), 
$$
therefore we have equality between all these expressions:
\begin{theorem}
\label{thm_1}
Assume that $H(\nu|\mu)<\infty$, where $d\nu=(E[e^{-f}])^{-1}
e^{-f}d\mu$ and $f$ is a measurable function. Then the infimum 
\beaa
J_f^\star&=&\inf(J_f(\xi):\xi\in \DD_{2,0}(H))\\
&=&\inf\left(\int f\circ (I_W+\xi)
  d\mu+H((I_W+\xi)\mu|\mu):\,\xi\in \DD_{2,0}(H)\right)
\eeaa
is attained at the  vector field $\xi=\nabla\varphi$, where
$\varphi$ is the unique (up to an additive constant) Monge potential
such that $(I_W+\nabla\varphi)\mu=\nu$ and that the $L^2(\mu,H)$-norm of $\nabla\varphi$
is equal to the Wasserstein distance between $\nu$ and $\mu$:
\beaa
d_H^2(\mu,\nu)&=&\inf\left(\int_{W\times W}|x-y|_H^2d\ga(x,y):\,\ga\in
  \Sigma_1(\mu,\nu)\right)\\
&=&\int_W |\nabla\varphi|_H^2 d\mu
\eeaa
where $\Sigma_1(\mu,\nu)$ denotes the set of probability measures on
$W\times W$, whose first marginals are $\mu$ and the second ones are
$\nu$.
\end{theorem}

 \noindent
Note that if we could apply the variational principle above, namely,
by taking the derivative of the functional $J_f$ at the minimizing
vector field $\nabla\varphi$ in any admissible direction, we would
obtain the following relation:
$$
\delta((I_H+\nabla^2\varphi)^{-1}-I_H)=\nabla\varphi+\nabla f\circ
(I_W+\nabla\varphi), 
$$
where $\delta$ denotes the Gaussian divergence, i.e., the adjoint of
the derivative $\nabla$ w.r.t. the Gaussian measure $\mu$ and this
equation implies 
Sobolev regularity of $\varphi$. A similar
method can be used for the dual Monge potential $\psi$ also. We shall
realize this program in the sequel beginning from the finite
dimensions and passing to the infinite dimensional case by a limiting
argument. This limit procedure requires more general
approximation-stability results about the convergence of Monge
potentials corresponding to convergent sequences of target measures
than those studied in \cite{fandu2}, they are delicate  (cf., Lemmas
\ref{OU-lemma} and \ref{min-lemma}) and they are of independent interest.

 Let us note that this method is applicable in other
situations than the Gaussian case as one can see already in the
case of dual potential. 

We make a last remark: this work is
devoted to the creation of a variational calculus by parametrizing the
formula \ref{1st-variation} with the vector fields which are
derivatives of scalar functionals. In another work, which has already
appeared, \cite{ASU-4}, we have parametrized the same formula with adapted vector
fields to obtain totally different results, like the existence, 
uniqueness and non-existence results of stochastic differential
equations with past depending drift coefficients.

\section{\bf{Preliminaries}}
Let $W$ be a separable Fr\'echet space equipped with  a Gaussian
measure $\mu$ of zero mean whose support is the whole
space\footnote{The reader may assume that $W=C(\R_+,\R^d)$, $d\geq 1$
  or $W=\R^{\N}$.}. The
corresponding Cameron-Martin space is denoted by $H$. Recall that the
injection $H\hookrightarrow W$ is compact and its adjoint is the
natural injection $W^\star\hookrightarrow H^\star\subset
L^2(\mu)$. The triple $(W,\mu,H)$ is called 
an abstract Wiener space. Recall that $W=H$ if and only if $W$ is
finite dimensional. A subspace $F$ of $H$ is called regular if the
corresponding orthogonal projection 
has a continuous extension to $W$, denoted again  by the same letter.
It is well-known that there exists an increasing sequence of regular
subspaces $(F_n,n\geq 1)$, called total,  such that $\cup_nF_n$ is
dense in $H$ and in $W$. Let $V_n$  be the
$\sigma$-algebra generated by $\pi_{F_n}$, then  for any  $f\in
L^p(\mu)$, the martingale  sequence 
$(E[f|V_n],n\geq 1)$ converges to $f$ (strongly if 
$p<\infty$) in $L^p(\mu)$. Observe that the function
$f_n=E[f|V_n]$ can be identified with a function on the
finite dimensional abstract Wiener space $(F_n,\mu_n,F_n)$, where
$\mu_n=\pi_n\mu$. 

Since the translations of $\mu$ with the elements of $H$ induce measures
equivalent to $\mu$, the G\^ateaux  derivative in $H$ direction of the
random variables is a closable operator on $L^p(\mu)$-spaces and  this
closure will be denoted by $\nabla$ cf.,  for example
\cite{ASU}. The corresponding Sobolev spaces 
(the equivalence classes) of the  real random variables 
will be denoted as $\DD_{p,k}$, where $k\in \NN$ is the order of
differentiability and $p>1$ is the order of integrability. If the
random variables are with values in some separable Hilbert space, say
$\Phi$, then we shall define similarly the corresponding Sobolev
spaces and they are denoted as $\DD_{p,k}(\Phi)$, $p>1,\,k\in
\NN$. Since $\nabla:\DD_{p,k}\to\DD_{p,k-1}(H)$ is a continuous and
linear operator its adjoint is a well-defined operator which we
represent by $\delta$. In the case of classical Wiener space, i.e.,
when $W=C(\reals_+,\reals^d)$, then $\delta$ coincides with the It\^o
integral of the Lebesgue density of the adapted elements of
$\DD_{p,k}(H)$ (cf.\cite{ASU}). 

For any $t\geq 0$ and measurable $f:W\to \reals_+$, we note by
$$
P_tf(x)=\int_Wf\left(e^{-t}x+\sqrt{1-e^{-2t}}y\right)\mu(dy)\,,
$$
it is well-known that $(P_t,t\in \reals_+)$ is a hypercontractive
semigroup on $L^p(\mu),p>1$,  which is called the Ornstein-Uhlenbeck
semigroup (cf.\cite{ASU}). Its infinitesimal generator is denoted
by $-\calL$ and we call $\calL$ the Ornstein-Uhlenbeck operator
(sometimes called the number operator by the physicists). Due to the
Meyer inequalities (cf., for instance \cite{ASU}), the
norms defined by 
\begin{equation}
\label{norm}
\|\varphi\|_{p,k}=\|(I+\calL)^{k/2}\varphi\|_{L^p(\mu)}
\end{equation}
are equivalent to the norms defined by the iterates of the  Sobolev
derivative $\nabla$. This observation permits us to identify the duals
of the space $\DD_{p,k}(\Phi);p>1,\,k\in\NN$ by $\DD_{q,-k}(\Phi')$,
with $q^{-1}=1-p^{-1}$, 
where the latter  space is defined by replacing $k$ in (\ref{norm}) by
$-k$, this gives us the distribution spaces on the Wiener space $W$
(in fact we can take as $k$ any real number). An easy calculation 
shows that, formally, $\delta\circ \nabla=\calL$, and this permits us
to extend the  divergence and the derivative  operators to the
distributions as linear,  continuous operators. In fact
$\delta:\DD_{q,k}(H\otimes \Phi)\to \DD_{q,k-1}(\Phi)$ and 
$\nabla:\DD_{q,k}(\Phi)\to\DD_{q,k-1}(H\otimes \Phi)$ continuously, for
any $q>1$ and $k\in \reals$, where $H\otimes \Phi$ denotes the
completed Hilbert-Schmidt tensor product (cf., for instance
\cite{ASU}). 
\\

\noindent
We shall use the following results about the divergence operator whose
proof in the regular case is straightforward, for extensions with less
regularity hypothesis, we refer the reader to Theorem B.6.4 in\cite{BOOK}.
\begin{lemma}
\label{div_property}
\begin{enumerate}
\item Assume that $\xi:W\to H$ is a smooth function and let $T=I_W+u$, with
$u\in \DD_{p,1}(H)$ for some $p>1$. Suppose that $T\mu\ll \mu$, then
the following identity holds $\mu$-a.s.:
$$
(\delta\xi)\circ T=\delta(\xi\circ T)+(\xi\circ
T,u)_H+\trace((\nabla\xi)\circ T.\nabla u)\,.
$$
\item The second moment of a divergence w.r.t. the Gauss measure is
  given by
$$
E[(\delta(\xi))^2]=E[|\xi|_H^2]+E[\trace(\nabla\xi\cdot\nabla\xi)]\,,
$$
where $``\cdot"$ is the composition operation between two linear
(Hilbert-Schmidt) operators.
\end{enumerate}
\end{lemma}

The following assertion which has been proved by H. Sugita
(cf. \cite{HS} or \cite{ASU}) is useful : assume that $(Z_n,n\geq
1)\subset \DD'$ converges to $Z$ in $\DD'$, assume further that each
each $Z_n$ is a probability measure on $W$, then $Z$ is also a
probability and $(Z_n,n\geq 1)$ converges to $Z$ in the weak topology
of measures. In particular, a lower bounded distribution (in the sense
that there exists a constant $c\in \R$ such that $Z+c$ is a positive
distribution) is a (Radon) measure on $W$.

A  measurable  function
$f:W\to \reals\cup\{\infty\}$ is called $H$-convex (cf.\cite{fandu0}) if 
$$
h\to f(x+h)
$$
is convex $\mu$-almost surely, i.e., if for any $h,k\in H$, $s,t\in
[0,1],\,s+t=1$, we have 
$$
f(x+sh+tk)\leq sf(x+h)+tf(x+k)\,,
$$
almost surely, where the negligeable set on which this inequality
fails may depend on the choice of $s,h $ and of $k$. We can rephrase
this property by saying that $h\to(x\to f(x+h))$ is an $L^0(\mu)$-valued
convex function on $H$. 
$f$ is called $1$-convex if the map 
$$
h\to\left(x\to f(x+h)+\frac{1}{2}|h|_H^2\right)
$$
is convex on the Cameron-Martin space $H$ with values in
$L^0(\mu)$. Note that all these  notions are  compatible with the
$\mu$-equivalence classes of random variables thanks to the
Cameron-Martin theorem. It is proven in \cite{fandu0} that 
this definition  is equivalent  the following condition:
  Let $(\pi_n,n\geq 1)$ be a sequence of regular, finite dimensional,
  orthogonal projections of $H$,  increasing to the identity map
  $I_H$. Denote also  by $\pi_n$ its  continuous extension  to $W$ and
  define $\pi_n^\bot=I_W-\pi_n$. For $x\in W$, let $x_n=\pi_nx$ and
  $x_n^\bot=\pi_n^\bot x$.   Then $f$ is $1$-convex if and only if 
$$
x_n\to \frac{1}{2}|x_n|_H^2+f(x_n+x_n^\bot)
$$ 
is  $\pi_n^\bot\mu$-almost surely convex. 
We define similarly the notion of $H$-concave and $H$-log-concave
functions. In particular, one can prove that, for any $H$-log-concave
function $f$ on $W$, $P_tf$ and $E[f|V_n]$ are again $H$-log-concave
\cite{fandu0}.

\section{\bf{Variational calculations}}

\noindent
Assume for a while that $\varphi\in\DD_{2,1}$ is smooth; this can be
achived by replacing $f$ by its regularization defined as 
$$
e^{-f_n}=E[P_{1/n}e^{-f}|V_n]\,,
$$
where $(P_t,t\geq 0)$ is the Ornstein-Uhlenbeck semi-group, $V_n$
is the sigma-algebra generated by $\{\delta e_1,\ldots,\delta e_n\}$
and $(e_n,n\geq 1)$ is a complete, orthonormal basis of $H$. Since
$J_f^\star=J(\nabla\varphi)$, if we take the Gateau derivative of $J$ at
$\nabla\varphi$, it should give zero: Let $L=(E[^{-f}])^{-1}e^{-f}$
and denote by $\La$ the Gaussian Jacobian of $T=I_W+\nabla\varphi$:
\begin{equation}
\label{lambda}
\La=\dett(I_H+\nabla^2\varphi)\exp\left(-\calL\varphi-\half|\nabla\varphi|_H^2\right)
\end{equation}
where $\calL$ is the Ornstein-Uhlenbeck operator $\calL=\delta\circ
\nabla$, $\dett$ denotes the modified Carleman-Fredholm determinant,
$\delta=\nabla^\star$ where the  adjoint is taken
w.r.t. the Wiener measure $\mu$, c.f. \cite{BOOK}. It follows from the
change of variables formula, c.f.\cite{BOOK}, that
$L\circ(I_W+\nabla\varphi)\,\La=1$, hence
$$
H((I_W+\nabla\varphi)\mu|\mu)=E\left[\half|\nabla\varphi|_H^2-\log
  \dett(I_H+\nabla^2\varphi)\right]\,.
$$
In particular, thanks to the $1$-convexity of $\varphi$, if we replace
$\varphi$ by $t\varphi$, for small $t\in [0,1]$, the shift
$T_t=I_W+t\nabla\varphi$ becomes strongly monotone and it is the
solution of the Monge transportation problem for the measure
$\nu_t=T_t\mu$ (i.e., the image of $\mu$ under $T_t$). Let $f_t$ be
defined as
$$
L_t=\frac{d\nu_t}{d\mu}=ce^{-f_t}\,.
$$
 If $\xi\in\DD_{2,1}(H)$ such that  $\nabla\xi$ has
small $L^\infty$-norm as a Hilbert-Schmidt operator, then
$T_{t,\eps}=I_W+t\nabla\varphi+\eps\xi$ is a strongly monotone shift for
small $t,\eps>0$, hence it
is almost-surely invertible (cf. \cite{BOOK}, Corollary 6.4.2). Note
moreover that the shift $I_W+t\nabla \varphi$ is the unique solution of
another Monge problem, namely the one which corresponds to the measure
$ce^{-f_t}d\mu$. Here the multiplication with a small $t$ permits us
to have a sufficiently large set on which we calculate the Gateau
derivative while preserving the $1$-convexity of the corresponding
Monge potential, namely $t\varphi$. 
Using again the change of variables formula for $T_{t,\eps}$, we get
$$
H(T_{t,\eps}\mu|\mu)=E\left[\half|t\nabla\varphi+\eps\xi|_H^2-\log
  \dett(I_H+t\nabla^2\varphi+\eps\nabla\xi)\right]\,.
$$
Therefore
$$
J_t(t\nabla\varphi+\eps\xi)=E\left[f_t\circ T_{t,\eps}+\half\left|t\nabla\varphi+\eps\xi\right|_H^2-\log
  \dett\left(I_H+t\nabla^2\varphi+\eps\nabla\xi\right)\right]\,.
$$
Since $t\nabla\varphi$ minimizes the function $J_t$ between all the
absolutely continuous shifts, we should have
\bea
\label{eqn_1}\lefteqn{\frac{d}{d\eps}J_t(t\nabla\varphi+\eps\xi)|_{\eps=0}}\\
&=&E\left[(t\nabla\varphi,\xi)_H-\trace\left(((I+t\nabla^2\varphi)^{-1}-I)\cdot(\nabla\xi)\right)+(\nabla
f_t\circ(I_W+t\nabla\varphi),\xi)_H\right]\nonumber\\
&=& 0\nonumber
\eea
for any $\xi\in \DD_{2,1}(H)$ with
$\|\nabla\xi\|_{2}\in L^\infty(\mu)$. Since the set of vector fields
$$
\Theta=\{\xi\in \DD_{2,1}(H):\,\|\nabla\xi\|_2\in L^\infty(\mu)\}
$$ 
is dense in any $L^p(\mu,H)$, we have proved the following
\begin{theorem}
\label{basic-thm}
In the finite dimensional smooth case, the Monge potential $\varphi$
satisfies the following relation
$$
\nabla\varphi+\nabla f\circ
(I_W+\nabla\varphi)-\delta\left[(I_H+\nabla^2\varphi)^{-1}-I_H\right]=0
$$
almost surely, where $\delta$ denotes the Gaussian divergence
w.r.t. $\mu$, i.e., the adjoint of $\nabla$ w.r.t. $\mu$.
\end{theorem}
\nproof
In the equation (\ref{eqn_1}) we have a term with trace, we just
interpret it as a scalar product on the Hilbert-Schmidt operators on
the Cameron-Martin space and the claim follows, for the case
$t\varphi$, from the definition of $\delta$ as a mapping from
Hilbert-Schmidt-valued operators to the vector fields under this
scalar product. Hence we have the identity
\begin{equation}
\label{eqn_2}
t\nabla\varphi+\nabla f_t\circ
(I_W+t\nabla\varphi)-\delta\left[(I_H+t\nabla^2\varphi)^{-1}-I_H\right]=0\,.
\end{equation}
Since we have $\La_t ce^{-f_t\circ T_t}=1$ a.s., where
$\La_t=\dett(I_H+t\nabla^2\varphi)\exp\left(-t\calL\varphi-\half|t\nabla\varphi|_H^2\right)$
and $T_t=I_W+t\nabla\varphi$, $\lim_{t\to 1}\nabla f_t\circ T_t=\nabla
f\circ T$ in probability, where $T=T_1=I_W+\nabla\varphi$. The
justification of the other terms being trivial, the proof is completed.
\nqed
\begin{theorem}
\label{control}
Suppose that $f$ is a smooth, bounded function on $\R^d$,
denote the probability measure $e^{-f}d\beta$ 
by $\nu$, where
$\beta$ is the standard Gauss measure on $\R^d$. Assume that the
forward and backward potentials associated to the transport couple
$(\beta,\nu)$  denoted respectively 
$(\varphi,\psi)$  are smooth functions. Let
$T$ be defined as
$T=I_{\R^d}+\nabla\varphi$. We have then
the following control:
\begin{eqnarray}
E\left[\left\|\nabla^2\psi\circ T\right\|_2^2\right]
&+&\sum_i E\left[\left\|(I+\nabla^2\varphi)^{1/2}\nabla^3\psi\circ
T(I+\nabla^2\varphi)^{1/2}e_i\right\|_2^2\right]\label{main_control}\\
&=&E\left[|\nabla\varphi+\nabla f\circ T|^2\right]\,,\nonumber
\end{eqnarray}
where $(e_i,i\leq d)$ is an orthonormal basis of $\R^d$ and
$\|\cdot\|_2$ denotes the Hilbert-Schmidt norm.
\end{theorem}
\nproof
Let $M=(I+\nabla^2\varphi)^{-1}-I$,
from Theorem \ref{basic-thm}, we have
$$
\delta(M)=\nabla\varphi+\nabla f\circ T
$$
and from Lemma \ref{div_property}
\beaa
\lefteqn{E[|\delta(M)|_{\R^d}^2]=\sum_{i=1}^d E[|\delta(Me_i)|^2]}\\
&&=\sum_{i=1}^d E\left[|Me_i|_{\R^d}^2+\trace((\nabla Me_i).(\nabla Me_i))\right]\,.
\eeaa
From the relation $(I+\nabla\psi)\circ T=I_{\R^d}$, we  obtain
$$
(I+\nabla^2\psi)\circ T=(I+\nabla^2\varphi)^{-1}\,,
$$
hence, differentiating both sides again we get
\beaa
\nabla((I+\nabla^2\psi)\circ T)&=&\nabla^3\psi\circ
T(I+\nabla^2\varphi)\\
&=&\nabla(I+\nabla^2\varphi)^{-1}\\
&=&-(I+\nabla^2\varphi)^{-1}\nabla^3\varphi
(I+\nabla^2\varphi)^{-1}\,,
\eeaa
this gives the identity
\begin{equation}
\label{eq_imp}
-(I+\nabla^2\varphi)^{-1}\nabla^3\varphi
(I+\nabla^2\varphi)^{-1}=\nabla^3\psi\circ T(I+\nabla^2\varphi)
\end{equation}
almost surely. 
Using the equality (\ref{eq_imp}) in the calculation of $\nabla M
e_i$, we get
\beaa
\nabla (Me_i)&=&-(I+\nabla^2\varphi)^{-1}\nabla^3\varphi
(I+\nabla^2\varphi)^{-1}e_i\\
&=&(\nabla_{e_i}\nabla^2\psi)\circ T(I+\nabla^2\varphi)\,.
\eeaa
We finally get
\beaa
\lefteqn{\trace((\nabla Me_i).(\nabla Me_i))}\\
&=&\trace\left((\nabla_{e_i}\nabla^2\psi)\circ T\,(I+\nabla^2\varphi)\cdot
  (\nabla_{e_i}\nabla^2\psi)\circ T \,(I+\nabla^2\varphi)\right)\\
&=&\trace\left(\left[(I+\nabla^2\varphi)^{1/2}(\nabla_{e_i}\nabla^2\psi)\circ T\,
(I+\nabla^2\varphi)^{1/2}\right]\left[(I+\nabla^2\varphi)^{1/2}(\nabla_{e_i}\nabla^2\psi)\circ T\,
(I+\nabla^2\varphi)^{1/2}\right]\right)\\
&=&\|(I+\nabla^2\varphi)^{1/2}(\nabla_{e_i}\nabla^2\psi)\circ T\, (I+\nabla^2\varphi)^{1/2}\|_2^2\,,
\eeaa
and the proof follows.
\nqed

\section{\bf{Approximations of the Monge Potentials}}
\label{approx_section}
The results  given in this section are indispensable to study the
stability and the approximation  results of the forward and backward potentials in the
finite dimensional situations whenever the target measures are
approximated with more regular measures. There are some results in the
literature (cg. \cite{CHE,Vil}), but they do not cover the situations
that we shall encounter.
\begin{lemma}
\label{OU-lemma}
Let $\beta$ be the standard Gaussian measure on $\R^d$,
$f\in\DD_{2,1}$ s.t. 
$$
\int_{\R^d}|\nabla f|^2e^{-f}d\beta<\infty\,.
$$
Let $(\varphi,\psi)$ be the Monge potentials associated to the
Monge-Kantorovitch problem  $\Sigma(\beta,\nu)$, where $d\nu=c
e^{-f}d\beta$. Define $f_n$ as to be $Q_{1/n}e^{-f}=e^{-f_n}$, where
$(Q_t,\,t\geq 0)$ denotes the Ornstein-Uhlenbeck semigroup on
$\R^d$. Let $(\varphi_n,\psi_n)$ be the Monge potentials corresponding
to Monge-Kantorovich problem  $\Sigma(\beta,\nu_n)$, where
$d\nu_n=ce^{-f_n}d\beta$. Then $(\varphi_n,n\geq 1)$ converges to
$\varphi$ in $\DD_{2,1}$, $(Q_{1/n}\psi_n,n\geq 1)$ converges to
$\psi$ in $L^1(\nu)$ and $(Q_{1/n}\nabla\psi_n,n\geq 1)$ converges to
$\nabla\psi$ in $L^2(\nu,\R^d)$
\end{lemma}
\nproof
In the sequel we replace $\varphi_n$ by $\varphi_n-E_\beta[\varphi_n]$
and $Q_{1/n}\psi_n$ by $Q_{1/n}\psi_n+E_\beta[\varphi_n]$ to avoid the
ambiguities about the constants.
Let $\ga_n,\ga$ be the unique solutions of Monge-Kantorovitch problems
for $(\beta,\nu_n)$ and $(\beta,\nu)$ respectively. From
Brenier's theorem (cf.\cite{BRE})
\begin{equation}
\label{Bre_1}
F_n(x,y)=\varphi_n(x)+\psi_n(y)+\half|x-y|^2=0 \,\,\,\ga_n-a.s.
\end{equation}
and $F_n(x,y)\geq 0$ for any $(x,y)\in \R^d\times\R^d$. Similarly
\begin{equation}
\label{Bre_2}
F(x,y)=\varphi(x)+\psi(y)+\half|x-y|^2=0 \,\,\,\ga-a.s.,
\end{equation}
and $F(x,y)\geq 0$ for any $(x,y)\in \R^d\times\R^d$.  Let us denote
by $p(x)=|x|^2$ (the Euclidean norm), then for any $\la<1/2$, we have 
$\int\,\exp[\la Q_{1/n} p(x)]d\beta(x)\leq \int\,\exp[\la p(x)]d\beta
(x)<\infty$ due to invariance of $\beta$ w.r.t. the Ornstein-Uhlenbeck
semigroup $(Q_t,t\geq 0)$ and due to the Jensen inequality. From the
definition of $\nu_n$
\begin{equation}
\label{variance}
\int |y|^2d\nu_n(y)=\int |y|^2 Q_{1/n}(e^{-f})d\beta\,.
\end{equation}
From the Young and Jensen inequalities, we obtain
\begin{equation}
\label{young-ineq}
|y|^2 Q_{1/n}(e^{-f})\leq\ e^{\eps|y|^2}+\frac{1}{\eps}
Q_{1/n}(e^{-f})\log Q_{1/n}(e^{-f}), 
\end{equation}
finally, again from the Jensen inequalities, it follows that the
sequence of
integrands at the right of the equality (\ref{variance}) is uniformly
integrable. Consequently it holds true that
$$
\lim_{n\to\infty}\int_{\R^d}|y|^2d\nu_n(y)=\int_{\R^d}|y|^2d\nu(y)
$$
and this implies that (cf. \cite{B-F}, Lemma 8.3)
\beaa
\lim_nE_\beta[|\nabla\varphi_n|^2]&=&\lim_nd_2(\beta,\nu_n)^2\\
&=&d_2(\beta,\nu)^2\\
&=&E_\beta[|\nabla\varphi|^2]\,,
\eeaa
where $d_2$ denotes the second order Wasserstein distance on the
probability measures on $\R^d$. These relations imply that
$(\varphi_n,n\geq 1)$ is bounded in $L^2(\ga)$. Moreover, taking into
account the relation $\nabla Q_t g=e^{-t} Q_t\nabla g$ for any smooth
$g$, we have
$$
E_{\nu_n}[|\nabla\psi_n|^2]=E_\beta[|\nabla\psi_n|^2Q_{1/n}e^{-f}]\geq
E_\beta[|\nabla Q_{1/n}\psi_n|^2e^{-f}]\,.
$$
By the boundedness of $(\varphi_n,n\geq 1)$ in $L^2(\beta)$ there
exists $a'\in L^2(\beta)$ such that  $(\varphi_n,n\geq 1)$ converges
weakly to $a'$ (upto a subsequence) in
$L^2(\beta)$, hence also in $L^2(\ga)$.  Since $F_n\geq 0$ everywhere,
by applying $Q_{1/n}$ in $y$-variable,  we get
$$
Q_{1/n} F_n(x,\cdot)(y)=\varphi_n(x)+Q_{1/n}\psi_n(y)+\half Q_{1/n}(|x-\cdot|^2)(y)\geq 0
$$
for any $(x,y)\in \R^d\times \R^d$. Integrating this inequality
w.r.t. the measure $\ga$ and taking the limit, we get
\beaa
\lefteqn{\lim_n\int\left(\varphi_n(x)+Q_{1/n}\psi_n (y)+\half
  Q_{1/n}(|x-\cdot|^2)(y)\right)d\ga}\\
&=&\lim_n\left(\int\varphi_nd\beta+\int\psi_nd\nu_n+\half
  Q_{1/n}(|x-\cdot|^2)(y) d\ga\right)\,.
\eeaa
Let us calculate  the third term:
\beaa
\int Q_{1/n}(|x-\cdot|^2)(y)
d\ga&=&\int\int|x-e^{-1/n}y-\sqrt{1-e^{-2/n}}z|^2d\beta(z)d\ga(x,y)\\
&=&\int\int\left(|x-e^{-1/n}y|^2+(1-e^{-2/n})|z|^2\right)d\beta(z)d\ga(x,y)
\eeaa
As 
$$
\lim_{n \to \infty}\int (1-e^{-1/n})|z|^2\beta(dz)=0
$$
and as 
\beaa
\int |y|^2d\ga&=&\int |y|^2d\nu=\int |y|^2 e^{-f} d\beta(y)\\
&\leq&\int e^{\eps|y|^2}d\beta (y)+\frac{1}{\eps} H(\nu|\beta)\,,
\eeaa
 we obtain
$$
\lim_n \int Q_{1/n}(|x-\cdot|^2)(y) d\ga=\int |x-y|^2 d\ga(x,y)\,.
$$
Let us note for later use that the inequality (\ref{young-ineq})
combined with the triangle inequality implies the $\ga$-uniform
integrability of the sequence $(Q_{1/n}(|x-\cdot|^2)(y),n\geq 1)$ and
from the dominated convergence theorem, we see that
$(Q_{1/n}(|x-\cdot|^2)(y),n\geq 1)$ converges to $|x-y|^2$ in
$\ga$-probability (i.e., in $L^0(\ga)$) since 
$$
|x-e^{-1/n}y+\sqrt{1-e^{-2/n}}z|^2\leq 4(|x|^2+|y|^2)+2|z|^2\,.
$$
\newline
By the Poincar\'e inequality, $(\varphi_n,n\geq 1)$ is bounded in $\DD_{2,1}$ w.r.t. the
Gaussian measure $\beta$, hence also it is bounded in $L^2(\ga)$. Therefore, upto a subsequence, it converges weakly
to some  $a'\in L^2(\ga)$ (note that $a'\in\DD_{2,1}(\beta)$ also). Moreover, from the relation (\ref{Bre_1})
\beaa
\lim_n\int Q_{1/n}\psi_n(y)d\ga(x,y)&=&\lim\int \psi_n(y) d\nu_n(y)\\
&=&-\lim\int\varphi_nd\beta-\half\lim\int |x-y|^2d\ga_n(x,y)\\
&=&-\int a d\beta-\half \int |x-y|^2d\ga(x,y)
\eeaa 
where the last equality above follows from
$$
\lim_n\half d_2^2(\beta,\nu_n)=\lim_n\half\int
|x-y|^2d\ga_n(x,y)=\int |x-y|^2d\ga(x,y)
$$
which is a consequence of  Lemma 8.3 of \cite{B-F}. Consequently
$$
\lim_n \int Q_{1/n}F_n(x,\cdot)(y)d\ga(x,y)=0\,.
$$
Therefore the sequence $\left(\varphi_n(x)+Q_{1/n}\psi_n(y)+\half|x-y|^2,\,n\geq 1\right)$
converges to $0$ in the norm topology of $L^1(\ga)$, consequently $(Q_{1/n}\psi_n,\,n\geq
1)$ is also uniformly integrable in $L^1(\ga)$, therefore there exists
some $b'\in L^1(\nu)$ which is a weak adherent point of $(Q_{1/n}\psi_n,\,n\geq
1)$. 
Therefore 
$$
a'(x)+b'(y)+\half|x-y|^2=0
$$
$\ga$-a.s. Let $(\varphi_n',n\geq 1)$ and $(Q_{1/n}\psi_n',n\geq 1)$
be the convex combinations of the sequences $(\varphi_n)$ and
$(Q_{1/n}\psi_n)$ respectively, which converge strongly in
$L^2(\ga)$ and $L^1(\ga)$ respectively. Let $a(x)=\limsup_n\varphi_n'(x)$ and
$b(y)=\limsup Q_{1/n}\psi_n'(y)$. We have then
$$
a(x)+b(y)+\half|x-y|^2\geq 0
$$
for all $(x,y)\in \R^d\times \R^d$ and as $a=a'$ and $b=b'$
$\ga$-almost surely, we have
$$
a(x)+b(y)+\half|x-y|^2= 0
$$
$\ga$-almost surely. By the uniqueness of the solution of the
Monge-Kantorovitch problem we should have $a=\varphi$ and $b=\psi$
$\ga$-a.s. Consequently $(\varphi_n, n\geq 1)$ converges weakly to
$\varphi$ in $\DD_{2,1}^\beta$. As
$\lim_nE_\beta[|\nabla\varphi_n|^2]=E_\beta[|\nabla\varphi|^2]$ and as
$(\nabla\varphi_n,n\geq 1)$ converges to $\nabla\varphi$ weakly, we
deduce the norm convergence of $(\varphi_n, n\geq 1)$ to $\varphi$ in
$\DD_{2,1}^\beta$. As $(Q_{1/n}(|x-\cdot|^2)(y),n\geq 1)$ converges to
$|x-y|^2$ in $L^1(\ga)$ (due to its $\ga$-uniform integrability and
its convergence in $L^0(\ga)$), we conclude that 
 $(Q_{1/n}\psi_n,n\geq 1)$ converges to $\psi$ (strongly) in $L^1(\nu)$, moreover $\nabla$ is closable
on $L^p(\nu),\,p\geq 1$ and
$\lim_nE_\nu[|\nabla Q_{1/n}\psi_n|^2]=E_\nu[|\nabla\psi|_H^2]$ and
this completes the proof.
\nqed

\noindent
In the next section we shall need the following result which is a
corollary of Lemma \ref{OU-lemma} whose notations are used without
further explanation:
\begin{corollary}
\label{OU-cor}
Assume the hypothesis of Lemma \ref{OU-lemma} are valid and assume
that $(\nabla^2\psi_n\circ T_n,n\geq 1)$
converges weakly in $L^2(\beta,\R^d\otimes\R^d)$, then $\nabla^2\psi\in
L^2(\nu,\R^d\otimes \R^d)$, where $\nabla$ is the closure in
$L^2(\nu)$ of the derivative operator, and the weak limit is equal to
$\nabla^2\psi\circ T$, i.e., 
$$
w-\lim_n \nabla^2\psi_n\circ T_n=\nabla^2\psi\circ T\,.
$$
\end{corollary}
\nproof
As $S\circ T=I_{\R^d}$ $\beta$-a.s., where
$S=I+\nabla\psi,\,T=I+\nabla\varphi$, we can represent
$\lim_n\nabla^2\psi_n\circ T_n$ as $\xi\circ T$, where $\xi\in
L^2(\nu,\R^d\otimes\R^d)$. As $\lim_n Q_{1/n}\psi_n=\psi$ in
$L^1(\nu)$, $\lim_nQ_{1/n}\nabla\psi_n=\nabla\psi$
and as $\nabla$ is closable in $L^2(\nu)${\footnote{This follows from
    $\int|\nabla f|^2d\nu<\infty$.}},
$(\nabla^2Q_{1/n}\psi_n,n\geq 1)$ converges as a distribution, i.e.,
for any smooth function of compact support $\eta$, we have 
$$
\lim_n
E_\nu[\nabla^2Q_{1/n}\psi_n\cdot\eta]=E_\nu[\nabla^2\psi\cdot\eta]\,.
$$
On the other hand, letting $L_n=Q_{1/n}L$,
\beaa
E_\nu[\nabla^2Q_{1/n}\psi_n\cdot\eta]&=&E[\nabla^2Q_{1/n}\psi_n\circ
T\cdot\eta\circ T]\\
&=&e^{-2/n}E[Q_{1/n}\nabla^2\psi_n\cdot \eta L]\\
&=&e^{-2/n}E[\nabla^2\psi_n\cdot Q_{1/n}(\eta L)]\\
&=&e^{-2/n}E\left[\nabla^2\psi_n\cdot \frac{Q_{1/n}(\eta L)}{L_n}L_n\right]\\
&=&e^{-2/n}E\left[\nabla^2\psi_n\circ T_n\cdot \left(\frac{Q_{1/n}(\eta
  L)}{L_n}\right)\circ T_n\right]. 
\eeaa
As 
$$
\left|\frac{Q_{1/n}(\eta L)}{L_n}\right|\leq\| \eta\|_{\infty}
$$
we get
\beaa
\lim_nE_\nu[\nabla^2Q_{1/n}\psi_n\cdot\eta]&=&E_\nu[\nabla^2\psi\cdot\eta]\\
&=&E[\xi\circ T\cdot\eta\circ T]\\
&=&E[\xi\cdot\eta L]\,,
\eeaa
consequently $\nabla\psi$
belongs to the $L^2(\nu)$-extended domain of the Gateaux derivative
operator  and  $\xi=\nabla^2\psi$ $\nu$-almost surely.
\nqed

\begin{lemma}
\label{mult-lemma}
Let $\beta$ be the standard Gaussian measure on $\R^d$,
$L\in L^1(\beta)$ be a probability density such that
$$
\int_{\R^d} L\log L d\beta<\infty\,.
$$
Let $(\varphi,\psi)$ be the Monge potentials associated to the
Monge-Kantorovitch problem  $\Sigma(\beta,\nu)$, where $d\nu=L
d\beta$. Define  $L_n(y)=c_nL(y)\tilde{\theta}_n(y)=\theta_n(y)L(y)$ as another density, where
$\tilde{\theta}_n\in C_K^\infty(\R^d)$ is approximating the constant
$1$, $c_n$ is the normalization constant and $\theta_n=c_n\tilde{\theta}_n$. Let $(\varphi_n,\psi_n)$
be the Monge potentials corresponding
to Monge-Kantorovich problem with quadratic cost over
$\Sigma(\beta,\nu_n)$, where $d\nu_n=L_nd\nu$. Then $(\varphi_n,n\geq
1)$ converges to
$\varphi$ in $\DD_{2,1}$, $(\theta_n\psi_n,n\geq 1)$ converges to
$\psi$ in $L^1(\nu)$ and also $(\sqrt{\theta_n}\nabla\psi_n, n\geq 1)$
converges to $\nabla\psi$ in $L^2(\ga)$, in particular
$$
\lim_nE_\nu[\theta_n|\nabla\psi_n|^2]=E_\nu[|\nabla\psi|^2]\,.
$$
\end{lemma}
\nproof
The proof is  similar to the proof of Lemma
\ref{OU-lemma}.
Let $\ga$ and $\ga_n$ be the transport plans
corresponding to the Monge-Kantorovitch problems for $(\beta,\nu)$ and
$(\beta,\nu_n)$ respectively. We replace  $\varphi_n$
by $\varphi_n-E_\beta[\varphi_n]$ and $\psi_n$ with
$\psi_n+E_\beta[\varphi_n]$ to fix the ideas. We denote also by
$\theta_n(y)$ the function $c_n\tilde{\theta}_n(y)$. 
As in the Lemma \ref{OU-lemma}, we have 
\begin{equation}
\label{allpts}
F_n(x,y)=\varphi_n(x)+\psi_n(y)+\half|x-y|^2\geq 0 
\end{equation}
for any $x,y\in \R^d$ and 
$$
F_n(x,y)=\varphi_n(x)+\psi_n(y)+\half|x-y|^2=0 \,\,\,\ga_n-a.s.
$$
It follows from (\ref{allpts}) that
\beaa
0&\leq&\int \theta_n(y)F_n(x,y)d\ga(x,y)\\
&=&\int\theta_n(y)\varphi_n(x)d\ga(x,y)+\int\theta_n(y)\psi_n(y)d\ga(x,y)\\
&&+\int\half\theta_n(y)|x-y|^2d\ga(x,y)\\
&=&I_n+II_n+III_n
\eeaa
Let us observe the second term in more detail:
\beaa
II_n&=&\int\psi_n(y)\theta_n(y)d\ga(x,y)=\int\psi_n(y)\theta_n(y)d\nu(y)\\
&=&\int\psi_n(y)d\nu_n(y)=\int\psi_n(y)d\ga_n(x,y)\\
&=&-\half\int|x-y|^2d\ga_n(x,y)-\int\varphi_n d\beta\,.
\eeaa
Consequently
\beaa
\int
\theta_n(y)F_n(x,y)d\ga(x,y)&=&\int\theta_n(y)\varphi_n(x)d\ga(x,y)
-\half\int|x-y|^2d\ga_n(x,y)\\
&&-\int\varphi_n d\beta+\half\int\theta_n(y)|x-y|^2d\ga(x,y)\\
&=&\int\theta_n(y)\varphi_n(x)d\ga(x,y)-\int\varphi_n(x)d\ga(x,y)\\
&&-\half\int|x-y|^2d\ga_n(x,y)\\
&&+\half\int\theta_n(y)|x-y|^2d\ga(x,y)\\
&=&\int(\theta_n(y)-1)\varphi_n(x)d\ga(x,y)\\
&&-\half\int|x-y|^2d\ga_n(x,y)
+\half\int\theta_n(y)|x-y|^2d\ga(x,y)
\eeaa
Recall that $(\tilde{\theta}_n, n\geq 1)$ is a sequence of positive smooth
functions of compact support increasing to one, hence the sequence
$(c_n, n\geq 1)$ defined as $c_n^{-1}=\int \tilde{\theta}_n(y)
L(y)d\beta(y)$, increases to one. Therefore, for any $\eps>0$ there
existes some $n_\eps\in \N$ such that $1\leq c_n\leq 1+\eps$ for any
$n\geq  n_\eps$. Remark also that, from Poincar\'e and Talagrand 
inequalities
\beaa
\int \varphi_n^2 d\ga&=&\int \varphi_n^2 d\beta\leq\int
|\nabla\varphi_n|^2d\beta\\
&\leq&2\int L_n\log L_nd\beta=2\int(\theta _n L\log L+\theta_n L\log
\theta_n)d\beta\\
&\leq&\frac{2}{1-\eps}\int L\log L d\beta+\frac{2}{1-\eps}\,,
\eeaa
where we have used the fact that $\theta_n\leq \frac{1}{1-\eps}$ for
$n\geq n_\eps$. As this estimation is uniform w.r.t. $n\geq n_\eps$,
we conclude that 
\begin{equation}
\label{phi_n}
\lim_{n\to\infty}\int\varphi_n(x)(\theta_n(y)-1)d\ga(x,y)=0\,.
\end{equation}
It follows from the monotone convergence theorem that 
\begin{equation}
\label{dist}
\lim_{n\to\infty}\int\theta_n(y) |x-y|^2d\ga(x,y)=\int
|x-y|^2d\ga(x,y)\,.
\end{equation}
The Young inequality can be applied as
$$
|y|^2\theta_n(y)L\leq\frac{1}{1-\eps}
(e^{\eps|y|^2}+\frac{1}{\eps}L\log L)\,,
$$
for $n\geq n_\eps$ as explained above. Hence the sequence
$(|y|^2\theta_n(y)L(y),n\geq 1)$ is uniformly integrable w.r.t. the
measure $\beta$, therefore
\begin{equation}
\label{bi}
\lim_{n\to\infty}\int |y|^2 L_nd\beta=\lim_{n\to\infty}\int |y|^2
d\nu_n=\int|y|^2d\nu\,.
\end{equation}
The relation (\ref{bi}), combined with \cite{B-F} implies
\beaa
\lim_n\half\int|x-y|^2d\ga_n(x,y)&=&\lim_n \half d_2^2(\nu_n,\beta)\\
&=&\half d_2^2(\nu,\beta)=\half\int|x-y|^2d\ga(x,y)\,.
\eeaa
The relations (\ref{phi_n}), (\ref{dist}) and (\ref{bi}) imply that
$$
\lim_n\int\theta_n(y) F_n(x,y)d\ga(x,y)=0\,.
$$
As $(\theta_n, n\geq n_\eps))$ converges to one and non-negative,
there exists a subsequence $(F_{n_k},k\geq 1)$ which converges to zero
$\ga$-a.s. Moreover $(\theta_n(y)\varphi_n(x),n\geq 1)$ is
$\ga$-uniformly integrable, hence $(\theta_n\psi_n, n\geq 1)$ is also
$\ga$-uniformly integrable. Consequently, upto a subsequence, the sequences
$(\theta_n(y)\varphi_n(x),n\geq 1)$ and $(\theta_n\psi_n, n\geq 1)$
converge weakly in $L^1(\ga)$ respectively to $a'$ and
$\psi'$. Moreover $(\varphi_n, n\geq 1)$ is bounded in the Sobolev
space $\DD_{2,1}^\beta$, hence it has a subsequence which converges to
some $\varphi'$ weakly in $\DD_{2,1}^\beta$, hence also weakly in
$L^2(\ga)$. For any $h\in L^\infty(\ga)$, we have
$$
\int(\varphi_n\theta_n-\varphi')hd\ga=\int\varphi_n(\theta_n-1)hd\ga+\int(\varphi_n-\varphi')hd\ga\,,
$$
As $(\theta_n, n\geq 1)$ converges to one in all $L^p$-spaces, the
first terms at the right converges to zero, the second one also
converges to zero as $n$ tends to infinity, therefore $a'=\varphi'$
$\ga$-a.s.
Let us take convex combinations of these weakly convergent sequences
to obtain the strong convergence these combinations, which assured by
Mazur's Lemma:
\begin{itemize}
\item The sequence of convex combinations
  $(co(\theta_n\varphi_n),n\geq 1)$ converges to $\varphi'$ in
  $L^2(\ga)$,
\item the sequence of convex combinations
  $(co(\theta_n\psi_n),n\geq 1)$ converges to $\psi'$ in
  $L^1(\ga)$,
\item $(\theta_n(y)|x-y|^2, n\geq 1)$ converges to $|x-y|^2$ in
  $L^1(\ga)$.
\end{itemize}
Define $\tilde{\varphi}=\lim\sup_n(co(\theta_n\varphi_n))$ and
$\tilde{\psi}=\lim\sup_n(co(\theta_n\psi_n))$. We have then
$$
\tilde{\varphi}(x)+\tilde{\psi}(y)+\half|x-y|^2\geq 0
$$
for any $x,y\in \R^d$ and we have also that
$$
\tilde{\varphi}(x)+\tilde{\psi}(y)+\half|x-y|^2= 0
$$
$\ga$-almost surely. By the uniqueness of the solutions of the dual
Monge-Kantorovitch, we should have $\varphi=\tilde{\varphi}$ and
  $\psi=\tilde{\psi}$ $\ga$-almost surely. As $(\varphi_n,n\geq 1)$
  converges to $\varphi$ weakly in $L^2(\beta)$ and as 
$\lim_nE_\beta[|\nabla\varphi_n|^2]=E_\beta[|\nabla\varphi|^2]$,
$(\varphi_n, n\geq 1)$ converges to $\varphi$ strongly in
$\DD_{2,1}^\beta$. As $(\theta_nF_n,n\geq 1)$ converges to zero
strongly in $L^1(\ga)$, $(\theta_n\psi_n, n\geq 1)$ converges strongly
to $\psi$ in $L^1(\ga)$. Moreover, we have 
\beaa
\int_{\R^d\times \R^d}\theta_n|\nabla\psi_n|^2d\ga&=&\int_{\R^d}\theta_n|\nabla\psi_n|^2d\nu\\
&=&\int|\nabla\psi_n|^2d\nu_n=\int|\nabla\psi_n\circ T_n|^2d\beta\\
&=&\int|\nabla\varphi_n|^2d\beta\,,
\eeaa
where $T_n=I_{\R^d}+\nabla \varphi_n$ is the forward transport map, and we obtain at once
\begin{equation}
\label{limit}
\lim_n\int_{\R^d\times
  \R^d}\theta_n|\nabla\psi_n|^2d\ga=\int|\nabla\psi|^2d\ga\,.
\end{equation}
In particular, the relation (\ref{limit}) implies the $\ga$-uniform
integrability of the sequence $(\sqrt{\theta_n}\nabla\psi_n, n\geq
1)$. We shall obtain the strong convergence of the latter as soon as
we prove that its weak limit is equal to $\nabla\psi$. To see this,
let $\xi$ be a smooth vector field, then 
\begin{equation}
\label{cvg}
\int(\nabla\psi_n,\xi)\sqrt{\theta_n}d\ga=\int_{\{\theta_n<c\}}(\nabla\psi_n,\xi)\sqrt{\theta_n}d\ga+\int_{\{\theta_n\geq
  c\}}(\nabla\psi_n,\xi)\sqrt{\theta_n}d\ga\,.
\end{equation}
As $\theta_n\to 1$, by the uniform integrability of
$(\sqrt{\theta_n}\nabla\psi_n, n\geq 1)$ and the boundeness of $\xi$, we see that the
first integral at the right of the equality (\ref{cvg}) can be made
arbitrarily small for $c<1$ for any $n\geq n_c$, for some $n_c\in
\N$. The second integral at the right of the equality (\ref{cvg}) is
equal to
$$
\int_{\{\theta_n\circ T_n\geq c\}}-\frac{(\nabla\varphi_n,\xi\circ
  T_n)}{\sqrt{\theta_n\circ T_n}}d\beta\,.
$$
By the uniform integrability of $(L_n, n\geq 1)$, the sequence $(T_n,
n\geq 1)$ is equi-concentrated on compacta and consequently $\lim\theta_n\circ
T_n=1$ in $\beta$-probability, also $\lim_n\xi\circ T_n=\xi\circ T$ in
$\beta$-probability. The dominated convergence theorem and the
$L^2$-boundedness of $(\nabla\varphi_n, n\geq 1)$ imply at one that 
\beaa
\lim_{n\to\infty}\int_{\{\theta_n\circ T_n\geq c\}}-\frac{(\nabla\varphi_n,\xi\circ
  T_n)}{\sqrt{\theta_n\circ T_n}}d\beta&=&-\int(\nabla\varphi, \xi\circ  T)d\beta\\
&=&\int(\nabla\psi\circ T, \xi\circ T)d\beta\\
&=&\int(\nabla\psi, \xi)d\ga
\eeaa

\nqed

\noindent
We also have the following, where we use the same notations as in the
preceding lemmas:
\begin{lemma}
\label{min-lemma}
Assume that $d\nu=c e^{-f}d\beta$, where $f:\R^d\to \R\cup\{\infty\}$
is a measurable function with the property that $\nu(\R^d)=1$ and that
\begin{equation}
\label{hyp-1}
\int (|f|+|\nabla f|^2)d\nu<\infty.
\end{equation}
Define $f_n=f\wedge n$, $n\in \N$ and define
$d\nu_n=c_ne^{-f_n}d\beta$. Let $(\varphi_n,\psi_n)$  and
$(\varphi,\psi)$ be the Monge
potentials corresponding to the Monge-Kantorovitch problem over
$\Sigma(\beta,\nu_n)$  and over $\Sigma(\beta,\nu)$ respectively (with
quadratic cost). Then $(\varphi_n, n\geq 1)$ converges to $\varphi$ in
$\DD_{2,1}^\beta$ and $(\psi_n,n\geq 1)$ converges to $\psi$ in
$L^1(\nu)$ as well as $(\nabla\psi_n,n\geq 1)$ converges $\nabla\psi$
in $L^2(\nu)$.
\end{lemma}
\nproof
We apply the same conventions about the expectations of $\varphi_n$ as
in the preceding lemmas. Let $F_n(x,y)=\varphi_n(x)+\psi_n(y)+\half|x-y|^2$, also define
$F(x,y)=\varphi(x)+\psi_(y)+\half|x-y|^2$. 
Denote by $\ga$
and $\ga_n$  the Monge-Kantorovitch plans for $\Sigma(\beta,\nu)$ and
$\Sigma(\beta, \nu_n)$ respectively. From the properties of the
solutions of the dual problem, we know that $F_n(x,y)\geq 0$ for all
$(x,y)\in \R^d\times \R^d$ and $F_n=0$ $\ga_n$-a.e., similar
properties hold true for $\ga$ and $F$, i.e., $F\geq 0$ everywhere and
$F=0$ $\ga$-a.e. In the sequel, to simplify the notations, we shall
equate the normalizing constants $c$ and $c_n$ to unity. Let
$D_f=\{x\in \R^d:f(x)<\infty\}$, we have
\beaa
\int\psi_nd\ga&=&\int \psi_n e^{-f}d\beta=\int\psi_n
e^{-(f-f_n)}d\nu_n\\
&=&\int \psi_n
e^{-(f-f_n)}d\ga_n=\int\left[-\varphi_n(x)-\half|x-y|^2\right]e^{-(f-f_n)}d\ga_n\\
&=&-\int\varphi_n(x)
e^{-(f(y)-f_n(y))}d\ga_n-\int\half|x-y|^2e^{-(f-f_n)}d\ga_n\,.
\eeaa
Therefore, using the relation $(I\times T_n)\beta=\ga_n$ and inserting
the expression for $\int\psi_nd\ga$ that we have calculated just
above, we get
\beaa
\int F_nd\ga&=&\int\psi_nd\ga+\int\varphi_nd\beta+\half\int|x-y|^2d\ga(x,y)\\
&=&\int\varphi_n(x)\left(1-e^{-(f(y)-f_n(y))}\right)d\ga_n(x,y)-\half\int|\nabla\varphi_n|^2
(1-e^{-(f-f_n)\circ T_n})d\beta\\
&=&I_n+II_n\,,
\eeaa
written in the respective order. 
We have 
\beaa
I_n&=&\int_{D_f}\varphi_n(x)\left(1-e^{-(f(y)-f_n(y))}\right)d\ga_n(x,y)\\
&&+\int_{D_f^c}\varphi_n(x)d\ga_n(x,y)\,.
\eeaa
Since $(\varphi_n,n\geq 1)$ is bounded
in $L^2(\beta)$ and as 
\begin{equation}
\label{Pr-cvg}
\lim_{n\to\infty}1_{T_n^{-1}(D_f)}(f\circ T_n-f_n\circ T_n)=0
\end{equation}
in $L^0(\beta)$ and $\lim_n\beta(T_n^{-1}(D_f^c))=0$ as well,  we conclude that $\lim_n I_n=0$. To calculate
$\lim_nII_n$, we proceed similarly: we have
\beaa
\int |x-y|^2e^{-(f-f_n)(y)}d\ga_n&=&\int_{T_n^{-1}(D_f)}
|\nabla\varphi_n|^2e^{-(f-f_n)\circ T_n} d\beta\\
&=&\int
|\nabla\varphi_n|^2e^{-(f-f_n)\circ T_n} d\beta\,.
\eeaa
 From the relation \ref{Pr-cvg}, the exponential term converges to 1
 in probability (i.e., in $L^0(\beta)$) boundedly. It remains to show
 that the sequence $(|\nabla\varphi_n|^2,n\geq 1)$ is uniformly
 integrable (w.r.t. $\beta$). By the triangle inequality, it suffices
 to prove that the sequence of functions $(|T_n(x)|^2,n\geq 1)$ is
 $\beta$-uniformly integrable and to achive this, using the notation $L_n=e^{-f_n}$, we write:
\beaa
\int_{\{|T_n|>c\}}|T_n|^2d\beta&=&\int_{\{|x|>c\}}|x|^2
L_n d\beta\\
&\leq&\int_{\{|x|>c\}}
e^{\eps|x|^2}d\beta+\frac{1}{\eps}\int_{\{|x|>c\}}L_n\log L_nd\beta\,,
\eeaa
the first term at the last line converges to zero as $c\to \infty$
from Fernique's Lemma (cf. \cite{ASU,ASU-1}). It is also easy to see that $(L_n\log L_n,n\geq1)$ is $\beta$-uniformly
integrable, hence 
$$
\lim_{c\to\infty}\sup_n \int_{\{|T_n|>c\}}|T_n|^2d\beta=0\,.
$$
Finally we see that 
$$
\lim_{n\to\infty}\int F_n(x,y)d\ga(x,y)=\lim_n(I_n+II_n)=0\,.
$$
Since $F_n\geq 0$ for any $n\geq 1$, $\lim\int F_nd\ga=0$  implies the
convergence of $(F_n, n\geq 1)$ in $L^1(\ga)$. As $(\varphi_n,n\geq
1)$ is bounded in $L^2(\beta)$, we can choose a weak cluster point of
it, say $a'$, as $(\psi_n,n\geq 1)$ is $\ga$-uniformly integrable, it
has also a weak cluster point, say $b'$. This implies that 
$$
a'(x)+b'(y)+\half|x-y|^2=0
$$
$\ga$-a.s. Thanks to the Mazur's lemma, by taking convex combinations of these two sequences, we
can obtain strongly convergent sequences $(\varphi_n')$ and
$(\psi_n')$, converging 
to the same limits $a'$ and $b'$ respectively in $L^2(\ga)$ and in $L^1(\ga)$. Define $a$ and $b$ respectively as
$a=\lim\sup_n\varphi_n',b=\lim\sup_n\psi_n'$, then, $a=a'$, $b=b'$
$\ga$-a.s. and 
$$
a(x)+b(y)+\half|x-y|^2\geq 0
$$
for all $(x,y)\in \R^d\times \R^d$ and
$$
a(x)+b(y)+\half|x-y|^2= 0
$$
$\ga$-a.s. From the uniqueness (up to constants) of the Monge
potential functions, it follows that $a=\varphi$ and $b=\psi$
$\ga$-a.s. Moreover, the same construction works for any cluster
points of the sequences $(\varphi_n)$ and $(\psi_n)$, consequently
these two sequences have each a unique cluster point and this proves
the convergence in $L^2(\ga)$ of $(\varphi_n,n\geq 1)$ and the convergence in
$L^1(\ga)$ of $(\psi_n, n\geq 1)$. Since
$\lim_nE_\beta[|\nabla\varphi_n|^2]=E_\beta[|\nabla\varphi|^2]$,
$(\varphi_n,n\geq 1)$ converges to $\varphi$ in $\DD_{2,1}^\beta$. To
show the convergence of $(\nabla\psi_n)$, we have, as $\ga(D_f)=1$,  
\beaa
E_\ga[|\nabla\psi_n|^2]&=&E_\ga[1_{D_f}|\nabla\psi_n|^2]\\
&=&E_\nu[1_{D_f}|\nabla\psi_n|^2]\\
&=&E_{\nu_n}\left[1_{D_f}|\nabla\psi_n|^2e^{-(f-f_n)}\right]\\
&=&E_\beta\left[1_{D_f}\circ T_n|\nabla\varphi_n|^2e^{-(f-f_n)\circ
  T_n}\right]\,.
\eeaa
As $\nabla\varphi_n\to\nabla\varphi$ in $L^2(\beta)$ and as
$1_{D_f}\circ T_n\exp[-(f-f_n)\circ T_n]\to 1$ in probability as
$n\to \infty$, we conclude that 
\beaa
\lim_n E_\ga[|\nabla\psi_n|^2]&=&\lim_n E_\ga[1_{D_f}|\nabla\psi_n|^2]\\
&=&\lim_n E_\beta\left[1_{D_f}\circ T_n|\nabla\varphi_n|^2e^{-(f-f_n)\circ
  T_n}\right]\\
&=&\lim_nE_\beta[|\nabla\varphi_n|^2]=E_\beta[|\nabla\varphi|^2]\\
&=&E_\nu[|\nabla\psi|^2]=E_\ga[|\nabla\psi|^2]\,.
\eeaa
Assume that $\xi\in L^2(\nu)$ be any weak cluster point of the
sequence $(\nabla\psi_n,n\geq 1)$, i.e., let
$\xi=\lim_k\nabla\psi_{n_k}$. As the derivative operator is closed in
$L^2(\nu)$ due to the hypothesis (\ref{hyp-1}), we should have
$\xi=\nabla\psi$, where the latter is defined in the Sobolev sense.
Consequently $(\nabla\psi_n, n\geq 1)$ converges to $\nabla\psi$ in
$L^2(\nu)$.
\nqed

\noindent
The proof of the following is exactly as the proof of Corollary
\ref{OU-cor}, hence it is omitted:
\begin{corollary}
\label{min-cor}
Assume the hypothesis of Lemma \ref{min-lemma} are valid and assume
that $(\nabla^2\psi_n\circ T_n,n\geq 1)$
converges weakly in $L^2(\beta,\R^d\otimes\R^d)$, then $\nabla^2\psi\in
L^2(\nu,\R^d\otimes \R^d)$, where $\nabla$ is the closure in
$L^2(\nu)$ of the derivative operator, and we have
$$
w-\lim_n \nabla^2\psi_n\circ T_n=\nabla^2\psi\circ T\,.
$$
\end{corollary}

\noindent
Lemma \ref{min-lemma}, although it turns the density $L=e^{-f}$ into a
non-degenerate one, it spoils the convexity of $f$, in other words the
log-concavity of $L$. The following results heals this default:

\begin{lemma}
\label{add-lemma}
Assume that $L=e^{-f}$ and $f$ satisfy the properties common to the
lemmata of this section. Define $L_\eps$ as 
$$
L_\eps=\frac{\eps+L}{1+\eps}\,,
$$
where $\eps>0$. Let $d\nu_\eps=L_\eps d\beta$ and  $d\nu=L d\beta$ then 
$$
\sup_{\eps\in [0,1]}H(\nu_\eps|\beta)<\infty\,.
$$
Moreover the forward and backward Monge potentials
$(\varphi_\eps,\psi_\eps)$ associated to the quadratic transport
problem on $\Sigma(\beta,\nu_\eps)$ satisfy the following properties:
\begin{itemize}
\item $\lim_{\eps\to 0}\varphi_\eps=\varphi$ in the Gaussian Sobolev
  space $\DD_{2,1}$,
\item $\lim_{\eps\to 0}\psi_\eps=\psi$ in $L^1(\nu)$ and
  $\lim_{\eps\to 0}\nabla\psi_\eps=\nabla\psi$ in $L^2(\nu)$.
\end{itemize}
\end{lemma}
\nproof
Since 
$$
L_\eps=\frac{\eps}{1+\eps}+\frac{1}{1+\eps}L\,,
$$
from the convexity of the function $x\to x\log x$, we get
\beaa
H(\nu_\eps|\beta)&=&\int_{\R^d} L_\eps\log L_\eps d\beta\\
&\leq&\frac{1}{1+\eps}\int_{\R^d} L\log L d\beta\,,
\eeaa
hence 
\begin{equation}
\label{in_eps}
\sup_{\eps\in [0,1]}H(\nu_\eps|\beta)\leq H(\nu|\beta)\,.
\end{equation}
Let now $\ga_\eps$ and $\ga$ be optimal transport plans associated to
$\Sigma(\beta,\nu_\eps)$ and $\Sigma(\beta,\nu)$ respectively, and let
$F_\eps(x,y)=\varphi_\eps(x)+\psi_\eps(y)+\half|x-y|^2$, where, as
usual we replace $\varphi_\eps$ and $\psi_\eps$ by
$\varphi_\eps-E_\beta[\varphi_\eps]$ and
$\psi_\eps+E_\beta[\varphi_\eps]$ respectively. We have
$F_\eps(x,y)\geq 0$ for all $x,y\in \R^d$ and $F_\eps=0$
$\ga_\eps$-a.s. We claim that 
\begin{equation}
\label{cvence}
\lim_{\eps\to 0}\int_{\R^d\times \R^d}F_\eps(x,y)d\ga(x,y)=0\,.
\end{equation}
At first as in the proof of Lemma \ref{OU-lemma}, $\lim_\eps\int
|y|^2d\nu_\eps=\int |y|^2d\nu$, hence we get 
\beaa
\lim_{\eps\to 0}E[|\nabla\varphi_\eps|^2]&=&\lim_{\eps\to 0}d_W^2(\beta,\nu_\eps)\\
&=&d_W^2(\beta)=E[|\nabla\varphi|^2]\,.
\eeaa
To prove (\ref{cvence}) we write first
\beaa
\int F_\eps d\ga+\eps\int\psi_\eps d\beta&=&\int_{\R^d}
\psi_\eps(\eps+L)d\beta+\half\int_{\R^d\times \R^d} |x-y|^2d\ga(x,y)\\
&=&(1+\eps)\int\psi_\eps d\nu_\eps+\half\int |x-y|^2d\ga(x,y)\,.
\eeaa
As
\beaa
\lim_{\eps\to 0}(1+\eps)\int\psi_\eps d\nu_\eps&=&-\lim_{\eps\to
  0}\half\int |x-y|^2d\ga_\eps(x,y)\\
&=&-\lim_{\eps\to 0} d_W^2(\beta,\nu_\eps)=-d_W^2(\beta,\nu)\\
&=&-\half\int |x-y|^2d\ga(x,y)\,,
\eeaa
we conclude that 
\begin{equation}
\label{lim1}
\lim_{\eps\to 0}\left(\int F_\eps d\ga+\eps\int\psi_\eps d\beta\right)=0\,.
\end{equation}
Recall that, $F_\eps\geq 0$ everywhere in $\R^d\times \R^d$, hence it
is also positive on the diagonal, i.e.,
$\varphi_\eps(x)+\psi_\eps(x)\geq 0$ for any $x\in \R^d$, as
$\int\varphi_\eps d\beta=0$, we have $\int\psi_\eps d\beta\geq 0$,
consequently the relation (\ref{cvence}) is proved. Note that, due to
the inequality (\ref{in_eps}) and the Talagrand inequality, the set
$(\varphi_\eps,\eps\in [0,1])$ is bounded in $\DD_{2,1}$.  The rest of the
proof goes exactly as the proofs of other lemmata of this section,
i.e., use of Mazur's Lemma, etc., hence it is omitted.
\nqed

\noindent
We have also the obvious corollary of Lemma \ref{add-lemma}, whose
proof is very similar to that of Corollary \ref{OU-cor}:
\begin{corollary}
\label{add-cor}
Assume the hypothesis of Lemma \ref{add-lemma} are valid and assume,
with $\eps=1/n$, that $(\nabla^2\psi_n\circ T_n, n\geq 1)$
converges in $L^2(\beta,\R^d\otimes\R^d)$, then $\nabla^2\psi\in
L^2(\nu,\R^d\otimes \R^d)$, where $\nabla$ is the closure in
$L^2(\nu)$ of the derivative operator, and we have
$$
\lim_n \nabla^2\psi_n\circ T_n=\nabla^2\psi\circ T\,.
$$
\end{corollary}

\section{\bf{Sobolev regularity}}
\label{sob_reg}
To approximate the Monge potentials constructed on the Wiener space
$(W,H, \mu)$ (cf.\cite{fandu1,fandu2}), as explained in preliminaries, we choose any complete
orthonormal basis $(e_n,n\geq 1)$ in the Cameron-Martin space $H$ and
construct the increasing sequence of $\sigma$-algebras $(V_n,n\geq
1)$, where $V_n$ is the $\sigma$-algebra generated by the Gaussian
random variables $\{\delta e_1,\ldots,\delta e_n\}$. 

The following lemma  is placed in this section although it is an approximation
result, since it is more interesting in the infinite dimensional case
and it is proved in \cite{fandu2}, we recall it here
for the sake of completeness:

\begin{lemma}
\label{condex-lemma}
\begin{enumerate}
\item
Let $L\in L^1_+(\mu)$, with $E[L]=1$ and define
$d\nu=Ld\mu=e^{-f}d\mu$. Assume
that $E[L|\nabla f|^2]<\infty$ and denote by $(\varphi,\psi)$ forward and
backward potentials corresponding to the transport problem from $\mu$
to $\nu$. Let $\nu_n$ be defined as $d\nu_n=L_nd\mu$, $L_n=E[L|V_n]$
and denote by $(\varphi_n,\psi_n)$ the Monge potentials corresponding
to $(\mu,\nu_n)$. Then $(\varphi_n,n\geq 1)$ converges to $\varphi$ in
$\DD_{2,1}$, $(\psi_n,n \geq 1)$ converges to $\psi$ in $L^1(\nu)$ and
$(\nabla\psi_n,n\geq 1)$ converges to $\nabla\psi$ in $L^2(\nu,H)$. 
\item
In particular, letting $T_n$ and $T$ to be the forward  transport maps
defined as $T_n=I_W+\nabla\varphi_n$, $T=I_W+\nabla\varphi$,
 if $(\nabla^2\psi_n\circ T_n, n\geq 1)$ converges weakly in 
 $L^2(\mu,H\otimes H)$, then, $\nabla^2\psi\in L^2(\nu,
H\times H)$ and we have
$$
w-\lim_n\nabla^2\psi_n\circ T_n=\nabla^2\psi\circ T,
$$
$\mu$-almost surely.
\end{enumerate}
\end{lemma}
\nproof
Only the second part requires a proof but it is immediate due to the
closability of the $H$-Gateaux derivative w.r.t. the target measure
$\nu$.
\nqed

We then define
the smooth target measures approximating the original one in three
steps: 
\begin{enumerate}
\item $d\nu_n=L_nd\mu$ as $L_n=E[e^{-f}|V_n]=e^{-f_n}$, 
\item $d\nu_{n,k}=L_{n,k} d\mu=e^{-f_{n,k}}d\mu=c_{n,k}e^{-f_n\wedge k}d\mu$ 
\item
  $d\nu_{n,k,l}=L_{n,k,l}d\mu=e^{-f_{n,k,l}}d\mu=P_{1/l}(e^{-f_{n,k}})d\mu$,  
\end{enumerate}
where $(P_t,t\geq 0)$ denotes the Ornstein-Uhlenbeck semigroup on the
Wiener space, whose version on $\R^d$ was denoted by $(Q_t,t\geq 0)$
in the preceeding pages. Lets us note again the relation which will be
used without any warning in the sequel: if $G:W\to Z$ is a Wiener
function in $L^p(\mu,Z)$, where $Z$ is a separable Hilbert space, then
$E[G|V_n]$ can be represented as $g_n(\delta e_1,\ldots, \delta e_n)$,
where $g_n:\R^n\to Z$ is in $L^p(\mu_n,Z)$, and $\mu_n$ is the Gauss
measure on $\R^n$. In this situation we have
$$
E[P_tG|V_n]=(Q_t g_n)(\delta e_1,\ldots, \delta e_n)
$$
$\mu$-a.s., where $(Q_t,t\geq 0)$ is the Ornstein-Uhlenbeck defined on $\R^n$.
Consequently  $L_{n,k,l}$ can be written as $e^{-\tilde{f}_n(\delta
  e_1,\ldots,\delta e_n)}$ where $\tilde{f}_n:\R^n\to \R$ is a smooth, lower
and upper bounded function. Consequently the classical results of transport like
\cite{Caf,Vil} affirm the existence of smooth Monge potentials
corresponding to the Monge-Kantorovitch problem on the set
$\Sigma(\mu_{n},\tilde{\nu}_{n,k,l})$, where $\mu_n$ is the Gauss measure on
$\R^n$ and $d\tilde{\nu}_{n,k,l}=exp(-\tilde{f}_{nkl})d\mu_n$.{\bf{ In the sequel we shall pass from the scenario on the Wiener
space with $f_n$ and the  measure $\mu$ to the scenario on  $\R^n$
with $\tilde{f}_n(x_1,\ldots,x_n)$ and the measure $\beta=\mu_n$, where $n$
will denote the dimension,  without further explanation and we shall
omit in the notations the ``tilde'' symbol}} for redactional simplicity. We shall
denote by $(\varphi,\psi)$,
$(\varphi_n,\psi_n)$,$(\varphi_{n,k},\psi_{n,k})$ and
$(\varphi_{n,k,l},\psi_{n,k,l})${\footnote{In the sequel we shall often
  omit the commas between the lower indices for typographical reasons}}the forward and backward Monge
potentials corresponding to the couples of measures $(\mu,\nu)$,
$(\mu,\nu_n)$, $(\mu_n,\nu_{n,k})$ and $(\mu_n,\nu_{n,k,l})$
respectively. Recall that, from the series of lemmas of Section
\ref{approx_section}, we know that
$\lim_{l\to\infty}\varphi_{n,k,l}=\varphi_{n,k}$,
$\lim_{k\to\infty}\varphi_{n,k}=\varphi_{n}$ and
$\lim_{n\to\infty}\varphi_{n}=\varphi$ in the Sobolev spaces
$\DD_{2,1}$. For the dual or backward potentials the situation is more
involved. We begin with
\begin{lemma}
\label{f-approx}
Let $T_{nkl}$ be the transport map defined as
$T_{nkl}=I_W+\nabla\varphi_{nkl}$, then we have
$$
\lim_n\lim_k\lim_l\nabla f_{nkl}\circ T_{nkl}=\nabla f\circ T
$$
in $L^2(\mu,H)$.
\end{lemma}
\nproof
We have 
\beaa
E\left[|\nabla f_{nkl}|_H^2 e^{f_{nkl}}\right]&=&4 E\left[|\nabla
  e^{-f_{nkl}/2}|_H^2\right]\\
&=&4E\left[\left|\half\frac{\nabla P_{1/l}(e^{-f_{nk}})}{P_{1/l}(e^{-f_{nk}})^{1/2}}\right|_H^2\right]\\
&=&E\left[\left|\frac{P_{1/l}(e^{-f_{nk}})}{P_{1/l}(e^{-f_{nk}})^{1/2}}\right|_H^2\right]\\
&\leq&e^{-2/l}E\left[\frac{P_{1/l}(|\nabla
    e^{-f_{nk}}|^2)}{P_{1/l}(e^{-f_{nk}})}\right]\,,
\eeaa
as, due to the local character of the Sobolev derivative, $\nabla
f_{nk}=0$ $\mu$-a.s. on the set $\{f_{nk}=k\}=\{f_n\geq k\}$, the last
line of the above inequality can be upper bounded as follows
\beaa
&&e^{-2/l}E\left[P_{1/l}(|\nabla f_n1_{\{f_n\leq
        k\}}e^{-f_n\wedge k}|^2)
\left(P_{1/l}(e^{-f_{n}\wedge k})\right)^{-1}\right]\\
&\leq&e^{-2/l}E\left[P_{1/l}(|\nabla f_n1_{\{f_n\leq
        k\}}e^{-f_n}|^2)(P_{1/l}(e^{-f_{n}\wedge
        k}))^{-1}\right]\\
&\leq&E\left[P_{1/l}(|\nabla f_n|^2e^{-f_n}1_{\{f_n\leq
    k\}})\frac{P_{1/l}(e^{-f_n}1_{\{f_n\leq
      k\}})}{P_{1/l}(e^{-f_n\wedge k})}\right]\\
&\leq&e^{-2/l}E[P_{1/l}(|\nabla f_n|^2e^{-f_n})]\\
&=&e^{-2/l}4E[|\nabla e^{-f_n/2}|^2]\\
&=&e^{-2/l}E\left[\frac{|\nabla
    E[e^{-f}|V_n]|^2}{E[e^{-f}|V_n]}\right]\\
&\leq&e^{-2/l}E\left[\frac{|
    E[\nabla e^{-f}|V_n]|^2}{E[e^{-f}|V_n]}\right]\\
&\leq&e^{-2/l}E\left[\frac{1}{E[e^{-f}|V_n]}E[|\nabla
  f|^2e^{-f}|V_n]E[e^{-f}|V_n]\right]\\
&=&e^{-2/l}E[e^{-f}|\nabla f|^2]\,.
\eeaa
Therefore $(\nabla e^{- f_{nkl}/2},\,n\geq 1;k\geq 1;l\geq 1)$ is bounded
in $L^2(\mu,H)$. As $\lim_{n,k,l}f_{nkl}=f$ in
$L^0(\mu)${\footnote{Here the order of the limits is important}} and as
$E[|e^{-f_{nkl}/2}|^2]\to E[|e^{-f/2}|^2]$, $(\nabla
  e^{-f_{nkl}/2},\,n\geq 1;k\geq 1;l\geq 1)$ converges weakly to
  $\nabla e^{-f/2}$ in $L^2(\mu,H)$. Hence
\beaa
E[|\nabla e^{-f/2}|^2]&\leq&\lim\inf_{n,k,l}E[|\nabla
e^{-f_{nkl}/2}|^2]\\
&\leq& \lim\sup_{n,k,l}E[|\nabla e^{-f_{nkl}/2}|^2]\\
&\leq&E[|\nabla e^{-f/2}|^2]\,,
\eeaa
from the above calculations, consequently the strong convergence in
$\DD_{2,1}$ holds true. Finally we have 
\beaa
\lim_{n,k,l}E[|\nabla f_{nkl}\circ T_{nkl}|^2]&=&\lim_{n,k,l}E[|\nabla
f_{nkl}|^2e^{-f_{nkl}}]\\
&=&E[|\nabla f|^2e^{-f}]=E[|\nabla f\circ T|^2]\,,
\eeaa
moreover $\lim_{n,k,l}\nabla f_{nkl}\circ T_{nkl}=\nabla f\circ T$ in
$L^0(\mu)$ (i.e., in probability), hence the $L^2$-convergence
follows.
\nqed

\begin{theorem}
\label{psi_cvg_thm}
The sequence $(\nabla^2\psi_{nkl}\circ T_{nkl}:\,n\geq 1;k\geq 1;l\geq
1)$ converges to $\nabla^2\psi\circ T$ weakly in $L^2(\mu,H\otimes H)$ as
$l\to\infty$, then as $k\to\infty$ and then as $n\to \infty$. 
\end{theorem}
\nproof
The  proof follows from applications of Corollary
\ref{OU-cor}, Corollary \ref{min-cor} and Lemma \ref{condex-lemma}.
\nqed

\begin{lemma}
\label{u_i_lemma}
The set of functions $(\log L_{nkl}\circ T_{nkl}: n,k,l\in \N)$ is
uniformly integrable.
\end{lemma}
\nproof
We claim first that $(L_{nkl}\log L_{nkl}:\,n,k,l\in \N)$ is uniformly
integrable: in fact, it
follows from the Jensen inequality 
\beaa
L_{nkl}\log L_{nkl}&\leq&P_{1/l}E[L_k\log L_k|V_n]\\
&\leq&P_{1/l}E[L\log L 1_{\{f\leq k\}}+ke^{-k}1_{\{f>k\}}|V_n]\\
&\leq&P_{1/l}E[L\log L|V_n]+1\,,
\eeaa
as $L\log L\in L^1(\mu)$, the set $(P_{1/l}E[L\log L|V_n];\,n,l\in \N)$
is uniformly integrable, therefore $(L_{nkl}\log L_{nkl}:\,n,k,l\in
\N)$ is also uniformly integrable. To complete the proof it suffices
to see that
\beaa
E[\log L_{nkl}\circ T_{nkl}1_{\{\log L_{nkl}\circ
  T_{nkl}>c\}}]&=&E[L_{nkl}\log L_{nkl}1_{\{\log L_{nkl}>c\}}]\\
&=&E[L_{nkl}\log L_{nkl}1_{\{ L_{nkl}>e^c\}}]\\
&=&E[L_{nkl}\log L_{nkl}1_{\{ L_{nkl}\log L_{nkl}>ce^c\}}]\to 0
\eeaa 
as $c\to \infty$ uniformly w.r.t. $n,k,l\in \N$ by the uniform
integrability of $(L_{nkl}\log L_{nkl}:\,n,k,l\in \N)$.
\nqed

\begin{theorem}
\label{L-control}
The sequence $(\calL\psi_{nkl}\circ T_{nkl},n,k,l)$ converges to
$\calL\circ T$ in the sense of distributions. Moreover, $\calL\psi
\circ T$ is in fact an element of $L^1(\mu)$.
\end{theorem}
\nproof
We shall prove only the convergence of a subsequence, which is chosen
by a double diagonalization method. 
The  proof of convergence in the sense of distributions is
straightforward by duality. To show that $\calL\psi\circ T\in
L^1(\mu)$ or equivalently that $\calL\psi\in L^1(\nu)$ is more
delicate: as $(\nabla^2\psi_n\circ T_n, n\geq 1)$ converges weakly in
$L^2(\mu,H\otimes H)$, from Mazur's Theorem, we can form a sequence of
its convex combinations, 
denoted as $(\nabla^2\psi_n'\circ T'_n, n\geq 1)$, which
converges strongly in $L^2(\mu,H\otimes H)$, in particular, the
sequence 
$$
(\|\nabla^2\psi_n'\circ T'_n\|_2^2, n\geq 1)
$$
is uniformly integrable in $L^1(\mu)$.
Let us denote by
$(\calL\psi'_n\circ T_n',n\geq 1)$ the corresponding convex
combinations of $(\calL\psi_n\circ T_n,n\geq 1)$. We have 
\beaa
\calL\psi'_n\circ T_n'&\geq&\sum_i\la_i\left[-\log L_{n_i}\circ
T_{n_i}-\half|\nabla\psi_{n_i}\circ
T_{n_i}|_H^2-\half\|\nabla^2\psi_{n_i}\circ T_{n_i}\|_2^2\right]\\
&\geq&\sum_i\la_i\left[-\log L_{n_i}\circ
T_{n_i}-\half|\nabla\psi_{n_i}\circ
T_{n_i}|_H^2\right]-\half\|\nabla^2\psi'_n\circ T'_n\|_2^2\,,
\eeaa
where we have used the inequality $|\dett(I+A)|\leq \half \|A\|_2^2$,
for any Hilbert-Schmidt operator $A$.
From Lemma \ref{u_i_lemma}, the sequence $(\calL\psi'_n\circ T_n',n\geq 1)$ is lower
bounded by a uniformly integrable sequence. We also have 
\beaa
\calL\psi'_n\circ T_n'&=&\sum_i\la_i\calL \psi_{n_i}\circ T_{n_i}\\
&\leq&\sum_i \la_i\left[-\log L_{n_i}\circ
T_{n_i}+\log\dett(I+\nabla^2\psi_{n_i}\circ T_{n_i}) -\half|\nabla\psi_{n_i}\circ
T_{n_i}|_H^2\right]\\
&\leq&\sum_i \la_i\left[-\log L_{n_i}\circ
T_{n_i}-\half|\nabla\psi_{n_i}\circ
T_{n_i}|_H^2\right].
\eeaa
Hence, from again Lemma \ref{u_i_lemma}, the sequence $(\calL\psi'_n\circ T_n',n\geq 1)$ is also
upperbounded by a uniformly integrable sequence, therefore it is
itself uniformly integrable. Consequently, $(\calL\psi'_n\circ
T_n',n\geq 1)$ converges weakly in $L^1(\mu)$ and the limit is equal to $\calL\psi\circ T$.
\nqed

\remark 
\section{\bf{Disintegration Results}}
\label{dis_section}
\noindent
Let $(\varphi,\psi)$ be the forward and backward potentials
corresponding to the solution of the Monge-Kantorovitch problem
corresponding to $(\mu,\nu)$, let $\pi_n:H\to H_n$ be an orthogonal
projection from $H$ to a finite dimensional subspace of $H$, denoted
by $H_n$, we may assume the existence of a complete  orthonormal basis
of $H$, denoted $(e_k,k\geq 1)\subset W^\star$ with
$H_n={\mbox{span}}\{e_1,\ldots,e_n\}$. As each $e_i\in W^\star$, $\pi_n$
has a continuous extension to the whole space $W$ that we shall denote
again with the same notation. The following
result has been proved in \cite{fandu2} (Lemma 6.2), in much more general  case:

\begin{lemma}
\label{section-lemma}
$\ga\in\Sigma(\mu,\nu)$ be the optimal measure (i.e., the transport
plan). Let $\pi_n$ be defined as above and let
$\pi_n^\bot=I_W-\pi_n$ . Define  $p_n$ as  the 
projection from $W\times W$ onto $H_n$ with $p_n(x,y)=\pi_nx$ and let
$p_n^\bot(x,y)=\pi_n^\bot x$.  
Consider the Borel disintegration 
\beaa
\ga(\cdot)&=&\int_{H_n^\bot\times W}\ga(\,\cdot|x_n^\bot)\ga_n^\bot(dz_n^\bot)\\
&=&\int_{H_n^\bot}\ga(\,\cdot|x_n^\bot)\mu_n^\bot(dx_n^\bot)
\eeaa
along the projection  of $W\times W$ on $H_n^\bot$, where
$\mu_n^\bot$ is the measure
$\pi_n^\bot\mu$,  $\ga(\cdot\,|x_n^\bot)$ denotes the regular 
conditional probability $\ga(\cdot\,|p_n^\bot =x_n^\bot)$ and $\ga_n^\bot$
is the measure $p_n^\bot\ga$. Then,
$\mu_n^\bot$ and $\ga_n^\bot$-almost surely  
$\ga(\,\cdot|x_n^\bot)$ is optimal on $(x_n^\bot+H_n)\times W$, in the
sense that it realizes the quadratic Wasserstein distance between the
measures $P_1\ga(\cdot|x_n^\bot)$ and $P_2\ga(\cdot|x_n^\bot)$, where
$P_i$ are the projection maps defined on $(x_n^\bot+H_n)\times W$ by  $P_i(\xi_1,\xi_2)=\xi_i,\,i=1,2$
\end{lemma}

\noindent
The following result is just the Bayes' formula of classical
probability:
\begin{lemma}
\label{Bayes}
For any $g\in C_b(W)$, we have
\begin{equation}
\label{bayes_1}
E_\mu\left[g\circ
T|\pi_n^\bot=x_n^\bot\right]=\frac{E_\mu\left[gL|\pi_n^\bot=x_n^\bot\right]}{E_\mu[L|\pi_n^\bot=x_n^\bot]}=E_\nu\left[g|\pi_n^\bot=x_n^\bot\right]\,,
\end{equation}
in other words the image of the regular conditional probability
measure $\mu(\cdot|\pi_n^\bot=x_n^\bot)$ under the transport map $x_n\to
x_n+\pi_n\nabla\varphi(x_n^\bot+x_n)$ is equal to the regular conditional
probability measure $\nu(\cdot|\pi_n^\bot=x_n^\bot)$ and this latter
measure is absolutely continuous w.r.t. $\mu(\cdot|\pi_n^\bot=x_n^\bot)$ with the
corresponding Radon-Nikodym density given by 
$$
x_n\to \frac{L(x_n+x_n^\bot)}{E_\mu[L|\pi_n^\bot=x_n^\bot]}\,.
$$
\end{lemma}
\nproof
For typographical simplicity, we shall omit the lower index ``n'' in
the proof. Let now $g,h\in C_b(W)$, then from the very definition of
the conditional probability, we have
\beaa
\int g\circ T(x^\bot+x)\mu(dx|x^\bot) h(x^\bot)\mu^\bot(dx^\bot)&=&\int
E_\mu[g\circ T|\pi^\bot=x^\bot] h(\pi^\bot x)\mu(dx)\\
&=&\int g\circ T(x)h(\pi^\bot x)\mu(dx)\\
&=&\int (g\circ T)\, (h\circ\pi^\bot\circ S\circ T) d\mu\\
&=&\int g(x)h(\pi^\bot S(x))L(x)\mu(dx)\\
&=&\int E_\mu[g\,L|\pi^\bot S]h(\pi^\bot S) d\mu\\
&=&\int E_\mu[gL|\pi^\bot S]\frac{L}{E_\mu[L|\pi^\bot S]}h(\pi^\bot S) d\mu\\
&=&\int \frac{E_\mu[gL|\pi^\bot S]}{E_\mu[L|\pi^\bot S]}h(\pi^\bot S)d\nu\\
&=&\int\frac{E_\mu[gL|\pi^\bot=x^\bot]}{E_\mu[L|\pi^\bot=x^\bot]}h(x^\bot)d\mu^\bot
\eeaa
and the proof follows.
\nqed

\begin{lemma}
\label{fin-compo}
Let $t(w_n^\bot,\cdot):H_n\to H_n$ be defined as
$t(w_n^\bot,x)=x+\pi_n\nabla\varphi(w_n^\bot+x)$ and let
$s(w_n^\bot,\cdot):H_n\to H_n$ be defined as
$s(w_n^\bot,x)=x+\pi_n\nabla\psi(w_n^\bot+x)$. Then it holds that,
$\mu_n^\bot$-almost everywhere
\beaa
s(w_n^\bot,t(w_n^\bot,x))&=&x\,\,\mu(\cdot|\pi_n^\bot=w_n^\bot)-a.s.\\
t(w_n^\bot,s(w_n^\bot,x))&=&x\,\,\nu(\cdot|\pi_n^\bot=w_n^\bot)-a.s.
\eeaa
In particular, we have, for $\mu_n^\bot$-almost all $w_n^\bot$,
$$
\frac{L(w_n^\bot+x)}{E_\mu[L|\pi_n^\bot=w_n^\bot]}
=\dett(I_{H_n}+D^2\psi(w_n^\bot+x))\\
\exp[-\delta\pi_n\nabla\psi(w_n^\bot+x)-\half|\pi_n\nabla\psi(w_n^\bot+x)|^2]\,,
$$
$\mu(\cdot|\pi_n^\bot=w_n^\bot)$almost surely.
\end{lemma}
\nproof
Note that we have
\begin{equation}
\label{eqn_A}
\pi_n\nabla\psi\circ T(w)=-\pi_n\nabla\varphi(w)=-D\varphi(w_n^\bot+x)
\end{equation}
$\mu$ a.s., where $D$ is the derivative on $H_n$ which is  regarded as the
Euclidean space $\R^n$, $w_n^\bot=w-\pi_n(w)=w-x$. The left hand side
of the relation (\ref{eqn_A}) can also be written as 
$$
D\psi(w_n^\bot+t(w_n^\bot,x))\,,
$$
hence we get, for $\mu_n^\bot$-almost all $w_n^\bot$, 
$$
D\psi(w_n^\bot+t(w_n^\bot,x))=-D\varphi(w_n^\bot+x)
$$
$\mu(\cdot|\pi_n^\bot=w_n^\bot)$-almost surely. Consequently the
partial maps $t(w_n^\bot,\cdot)$ and $s(w_n^\bot,\cdot)$  are inverse
to each other on $H_n$, consequently the representation of the density
follows from finite dimensional results about the derivatives of the
convex functions and their Legendre transformations corresponding to
the Monge-Kantorovitch problem corresponding to the measures which are
absolutely continuous w.r.t. Lebesgue measure (cf.\cite{Mc,Vil}).
\nqed
\begin{remarkk}
From the definition of the conditional probability, we can represent
the Radon-Nikodym density of  Lemma \ref{fin-compo} as follows:
\bea
\label{key}
&&\frac{L(\pi_n^\bot(w)+x)}{E_\mu[L|\pi_n^\bot]}\\
&&=\dett(I+\pi_n\nabla^2\psi\pi_n(\pi_n^\bot(w)+x)\exp\left[-\delta\pi_n\psi(\pi_n^\bot(w)+x)-\half|\pi_n\nabla\psi(\pi_n^\bot(w)+x)|^2\right]\nonumber
\eea
$\mu_n^\bot\times \mu_n=\mu$-almost surely.
\end{remarkk}

\begin{lemma}
\label{ui_lemma}
The sequence $(\delta(\pi_n\nabla\psi),n\geq 1)$ is $\nu$-uniformly
integrable and it converges to $\calL\psi$ in $L^1(\nu)$. In
particular
$$
\lim_{n\to\infty}\frac{L(\pi_n^\bot(w)+\pi_n(w))}{E_\mu[L|\pi_n^\bot]}=\dett(I+\nabla^2\psi(w))
\exp\left[-\calL\psi(w)-\half|\nabla\psi(w)|^2\right]
$$
$\nu$-a.s.
\end{lemma}
\nproof
Let 
$$
\la_n(w_n,w_n^\bot)=\frac{L(w_n^\bot+w_n)}{E_\mu[L|w_n^\bot]}\,,
$$
from Martingale Convergence Theorem of Doob,
$(\la_n(\pi_n^\bot(w)+\pi_n(w)),n\geq 1)$ converges to $L$ on
$D=\{L>0\}$ $\mu$-a.s., or
$\nu$-a.s. Therefore $(\delta(\pi_n\nabla\psi),n\geq 1)$ converges
$\mu$-a.s. on $D$ (or $\nu$-a.s.). We have, by writing
$w=\pi_n(w)+\pi_n^\bot(w)$, from Lemma \ref{fin-compo}
\beaa
-\log\la_n(w)&=&-\log\dett(I+\pi_n\nabla^2\psi\pi_n)+\delta\pi_n\nabla\psi+\half|\pi_n\nabla\psi|_H^2\\
&=&-\log L+\log E[L|\pi_n^\bot]
\eeaa
$\mu$-a.s. on the set $D$. Consequently we can write
\begin{equation}
\label{expl_eq}
\delta\pi_n\nabla\psi=-\log L+\log
E[L|\pi_n^\bot]+\log\dett(I+\pi_n\nabla^2\psi\pi_n)-\half|\pi_n\nabla\psi|_H^2
\end{equation}
$\nu$-a.s. and $\mu$-a.s. on the set $D$. To show the uniform
integrability, we write, for any $A\in \calB(W)$
\beaa
E_\nu[1_A|\delta\pi_n\nabla\psi|]&\leq&E_\nu[-1_A\log
L]+E_\nu[1_A\,\log E[L|\pi_n^\bot]\\
&&+\half E_\nu[1_A\|\pi_n\nabla^2\psi\pi_n\|_2^2]+\half E_\nu[1_A|\pi_n\nabla\psi|_H^2]\\
&=&I+II_n+III_n+IV_n\,,
\eeaa
As $E_\mu[L\log L]<\infty$, we have $ I<\eps/4$, provided
$\mu(A)<\delta$. For $II_n$ we have, from Jensen inequality
\beaa
E_\nu\left[1_A\log E_\mu[L|\pi_n^\bot]\right]&=&
E_\nu\left[E_\nu[1_A|\pi_n^\bot] \log E_\mu[L|\pi_n^\bot]\right]\\
&=&E_\mu\left[E_\nu[1_A |\pi_n^\bot] E_\mu[L|\pi_n^\bot]\log E_\mu[L|\pi_n^\bot]\right]\\
&\leq&E_\mu\left[E_\nu[1_A|\pi_n^\bot] L\log L\right]
\eeaa
as $E_\nu[|\log L|]<\infty$, the sequence $(E_\nu[\log
L|\pi_n^\bot],\,n\geq 1)$ is  $\nu$-uniformly integrable, hence 
$\sup_nII_n\leq \eps/4$ for $\mu(A)<\delta$. As
$(\|\pi_n\nabla^2\psi\pi_n\|_2,n\geq 1)$ is a monotone, increasing
sequence, it follows from the Monotone Convergence Theorem that
$\sup_nIII_n<\eps/4$ for $\mu(A)<\delta$. The fourth term $IV_n$ is
trivial to control since $E_\nu[1_A|\pi_n\nabla\psi|_H^2]\leq
E_\nu[1_A|\nabla\psi|_H^2]$ for any $n\geq 1$. Consequently
$(\delta(\pi_n\nabla\psi),n\geq 1)$ is $\nu$-uniformly integrable,
hence $\lim_{n\to\infty} \delta(\pi_n\nabla\psi)=\calL\psi$ in
$L^1(\nu)$. 
\nqed

\section{\bf{Calculation of the Jacobians and Monge-Amp\`ere Equation}}
\noindent
We are now at a position to express the density $L$ in terms of the
backward potential for the original problem:
\begin{theorem}
\label{L_calcul}
The target density has the following representation:
$$
L=\dett(I+\nabla^2\psi)\exp\left[-\calL\psi-\half|\nabla\psi|_H^2\right]
$$
$\nu$-almost surely.
\end{theorem}
\nproof
Let us extract a sequence $(\Psi_n,n\geq 1)$ from $(\psi_{nkl})$ of
the form $\Psi_n=\psi_{n,k_n,l_{k_n}}$ applying twice the diagonal
sequence selection. By the uniform integrability results we can assume
that $w-\lim_n\nabla^2\Psi_n\circ T_n=\nabla^2\psi\circ
T,\,w-\lim_n\calL\Psi_n\circ T_n=\calL\psi\circ T$ weakly in
$L^2(\mu,H\otimes H)$ and weakly in $L^1(\mu)$ respectively, where
$T_n=I_W+\nabla\Phi_n$ and $\Phi_n=\varphi_{n,k_n,l_{k_n}}$. Note that
$\nabla\Psi_n\circ T_n$ converges strongly in $L^2(\mu,H)$. From
Mazur's theorem, we can construct sequence of convex combinations
\footnote{eventually taking twice convex combinations, the first combinations for assuring
the weak convergence of $(\calL\Psi_n\circ T_n,n\geq 1)$ in $L^1(\mu)$
and the  second ones to assure its strong convergence in $L^1(\mu)$ as
well as the strong convergence of $\nabla\Psi_n\circ T_n$ in $L^2(\mu,H\otimes H)$  }
$$
\left(\sum_{m_i\geq n}\la_i\nabla^2\Psi_{m_i}\circ T_{m_i},\,n\geq 1\right)
$$
which converges strongly to $\nabla^2\Psi\circ T$ in $L^2(\mu,H\otimes
H)$, hence
\beaa
&&\lefteqn{\lim_n\big[-\log\dett(I+\sum_{m_i\geq n}\la_i\nabla^2\Psi_{m_i}\circ T_{m_i})}\\
&&+\sum_{m_i\geq n}\la_i\calL\Psi_{m_i}\circ
  T_{m_i}
+\half\big|\sum_{m_i\geq n}\la_i\Psi_{m_i}\circ
  T_{m_i}\big|_H^2\big]\\
&=& -\log \dett(I+\nabla^2\psi\circ T)+\calL\psi\circ
T+\half|\nabla\psi\circ T|_H^2\\
&=&-\log\La(\psi)\circ T\,
\eeaa
in the weak topology. From the convexity of the $A\to
-\log\dett(I_H+A)$, it follows that 
$$
-\log\La(\psi)\circ T \leq -\log L\circ T\,,
$$
hence
$$
\La(\psi)\circ T\geq L\circ T
$$
$\mu$-almost surely or $\La(\psi)\geq L$ $\nu$-almost surely. To show
that they are equal $\nu$-a.e., it suffices to prove that 
\begin{equation}
\label{density-ineq}
E_\mu[1_D\La(\psi)]\leq 1\,,
\end{equation}
where $D$ is defined
as $D=\{w\in W:L(w)>0\}$. If the last claim were true we
would have 
$$
1=E[L1_D]\leq E[\La(\psi)1_D]\leq 1\,,
$$
which would imply  $L1_D=\La(\psi)1_D$ $\mu$-a.s., as
$\mu$ and $\nu$ are equivalent on $D$, the claim would have been
proved. Let us now prove (\ref{density-ineq}): From Lemma \ref{ui_lemma}, the Fatou Lemma and the Fubini Theorem, we have
\beaa
E_\mu[\La(\psi)]&\leq&\lim\inf_nE_\mu[\La(\pi_n^\bot+\pi_n))]\\
&=&\lim\inf_n\int\frac{\La(\psi(w_n^\bot+x))}{E_\mu[L|\pi_n^\bot=w_n^\bot]}d\mu_n(x)d\mu_n^\bot(w_n^\bot)=1
\eeaa
and the proof follows.
\nqed

Recall that we denote the forward potential of the original problem by
$\varphi$ which is a $1$-convex function, and it is an element of the
Sobolev space $\DD_{2,1}$. Its $1$-convexity implies that the operator
valued distribution $I_H+\nabla^2\varphi$ is positive, hence it is a
vector measure in the sense that, for any finite rank operator
$\kappa$ on $H$, $\trace[((I+\nabla^2\varphi)\kappa]$ is a signed measure
and if  $\kappa$ is also positive, then
$\trace((I+\nabla^2\varphi)\kappa)$ is a positive measure. Hence, due to
the Lebesgue decomposition theorem, we can
write it as the sum of an absolutely continuous measure
and a singular measure w.r.t. the Wiener measure $\mu$. The absolutely
continuous part is well-defined, denoted by $\nabla^2_a\varphi$ and it
satisfies
$$
[\trace((I+\nabla^2\varphi)\kappa)]_a=[\trace((I+\nabla^2_a\varphi)\kappa)]\,.
$$
Similarly, we can look at $\pi_n\nabla^2\varphi\pi_n$ as the second
order derivative in the direction of  the $n$-dimensional space
$H_n$, hence, if we denote $w=\pi_n^\bot(w)+\pi_n(w)=w_n^\bot+w_n$,
then the partial map $w_n^\bot\to
\pi_n\nabla^2\varphi\pi_n (w_n^\bot+\cdot)$ can be interpreted as a
measurable map with values in the space of measures on $H_n$, whose
absolutely continuous part w.r.t. Lebesgue measure is denoted as
$\pi_n\nabla^2\varphi\pi_n(w_n^\bot+\cdot)_a$ 
\footnote{Note that the lower index ``a'' depends also on the dimension
$n$ and we omit this dependence for typographical simplicity.}. With these notations we
can announce

\begin{lemma}
\label{ac_1_lemma}
The following relations hold true $\mu$-a.s.:
\begin{enumerate}
\item
$$
\lim_n (\pi_n\nabla^2\varphi\pi_n)_a=\nabla_a^2\varphi\,,
$$
\item
$$
\lim_n\left((\pi_n\nabla^2\varphi\pi_n)(\pi_n^\bot(w)+\cdot)\right)_a(\pi_n(w))=\nabla^2_a\varphi
$$
\end{enumerate}
\end{lemma}
\nproof
Let us recall that the inequalities that we use below are to be understood in the
sense of operators, i.e., we say that $A\geq B$, where $A$ and $B$ are
two symmetric, poistive operators on $H$ if $A-B$ is a positive
operator. As $(\pi_n,n\geq 1)$ increases to the identity operator of
$H$, the sequence $(\pi_n\nabla^2\varphi\pi_n,n\geq 1)$ increases to
$\nabla^2\varphi$ as a sequence of lower bounded operator-valued
measures, by the maximality of the absolutely continuous part of a
measure in Lebesgue decomposition theorem, it follows that
$\lim_n(\pi_n\nabla^2\varphi\pi_n)_a\leq \nabla^2_a\varphi$ almost
surely. We have also naturally $\nabla^2\varphi\geq
\nabla^2_a\varphi$, in particular, if $0\leq g\in \DD$ and if
$\gamma$ is a positive, symmetric, Hilbert-Schmidt or nuclear operator on $H$,
then 
$$
E[\trace(\nabla^2\varphi[g] \gamma)]\geq
E[\trace(\nabla^2_a\varphi\,\gamma)g]\,,
$$
and $\pi_n\nabla^2\varphi\pi_n\geq
\pi_n\nabla_a^2\varphi\pi_n$, by maximality we have also 
$$
( \pi_n\nabla^2\varphi\pi_n)_a\geq \pi_n\nabla_a^2\varphi\pi_n
$$
almost surely. Taking the limit of both sides we get
$$
\lim_n ( \pi_n\nabla^2\varphi\pi_n)_a\geq
\lim_n\pi_n\nabla_a^2\varphi\pi_n=\nabla^2_a\varphi\,,
$$which proves the first claim. To prove the second relation,  let
$w_n=\pi_nw,\,w_n^\bot=\pi_n^\bot w$. We have,
as explained above,
\beaa
(\pi_n\nabla^2\varphi\pi_n)_a(w_n^\bot+\cdot)&\leq&(\pi_n\nabla^2\varphi\pi_n)(w_n^\bot+\cdot)_a\\
&\leq&(\pi_n\nabla^2\varphi\pi_n)(w_n^\bot+\cdot)\leq
\nabla^2\varphi(w_n^\bot+\cdot)
\eeaa
$\mu_n$-a.s. for $\mu_n^\bot$-almost all $w_n^\bot$. Writing
$w_n=\pi_nw,\,w_n^\bot=\pi_n^\bot w$, via Fubini, we get
$$
\limsup_n(\pi_n\nabla^2\varphi\pi_n)(w_n^\bot+\cdot)_a\leq\nabla^2\varphi(w)
$$
$\mu$-a.s., and from the maximality of the absolutely continuous part,
we have also 
$$
\limsup_n(\pi_n\nabla^2\varphi\pi_n)(w_n^\bot+\cdot)_a\leq\nabla_a^2\varphi(w)\,,
$$
almost surely. Writing all the inequalities together we get
\beaa
\nabla_a^2\varphi&\geq&\limsup_n(\pi_n\nabla^2\varphi\pi_n)(w_n^\bot+\cdot)_a\\
&\geq&\limsup_n(\pi_n\nabla^2\varphi\pi_n)_a(w_n^\bot+\cdot)\\
&=&\limsup_n(\pi_n\nabla^2\varphi\pi_n)_a(w)\\
&=&\nabla^2_a\varphi
\eeaa
$\mu$-a.s.
\nqed

\noindent
As a consequence of the above result we have
\begin{theorem}
\label{a_regularity}
The second Alexandroff derivative of the forward transport potential
$\varphi$, which is denoted as $\nabla_a^2\varphi$,
is a map with values in the space of the Hilbert-Schmidt operators on
the Cameron-Martin space $H$. It satisfies the following identity:
\begin{equation}
\label{a1_identity}
(I_H+\nabla^2\psi\circ T)(I_H+\nabla^2_a\varphi)=I_H
\end{equation}
$\mu$-almost surely.
\end{theorem}
\nproof
From the finite dimensional results and with the notations explained
in Section \ref{dis_section} and in this section, cf., \cite{Mc,Vil}, we have, for
$\mu_n^\bot$-a.a. $w_n^\bot$ the relation
\begin{equation}
\label{a2_identity}
\left(I_{\R^n}+D_n^2\psi(w_n^\bot+\pi_nT(w_n^\bot+x))\right)(I_{\R^n}+(D_{n}^2\varphi(w_n^\bot+\cdot))_a(x))=I_{\R^n}\,,
\end{equation}
where, as we have mentioned before, $D_n$ is the derivation operator
on $n$-dimensional Euclidean space and the notation
$(D_{n}^2\varphi(w_n^\bot+\cdot))_a(x)$ means the second Alexandroff
derivative of the partial map $x\to \varphi(w_n^\bot+x)$ while
$w_n^\bot$ is kept fixed. From the equation (\ref{a2_identity}), we get
\begin{equation}
\label{a3_identity}
D_n^2\psi(w_n^\bot+\pi_nT(w_n^\bot+x))(I_{\R^n}+\pi_n\nabla^2\varphi\pi_n(w_n^\bot+\cdot)_a)(x)
=-\pi_n\nabla^2\varphi\pi_n(w_n^\bot+\cdot)_a\,.
\end{equation}
As 
$$
\lim_n
D_n^2\psi(\pi_n^\bot(w)+\pi_nT(\pi_n^\bot(w)+\pi_n(w))=\nabla^2\psi\circ
T(w)
$$
$\mu$-a.s., and 
$$
\lim_n(D_n^2\varphi(\pi_n^\bot+\cdot))_a(\pi_n(w)=\nabla_a^2\varphi(w)
$$
$\mu$-a.s., as proven in Lemma \ref{ac_1_lemma}, the proof follows
once we recall that $D_n$ is the restriction of $\nabla$ to the
$n$-dimensional subspaces of $H$.
\nqed

\begin{remarkk}
\label{cv_rem}
It follows from the last lines of the above proof, namely from the
relation (\ref{a3_identity}), that 
$$
\lim_n
\pi_n\nabla^2\varphi\pi_n(\pi_n^\bot(w+\cdot)_a)(\pi_n(w))=\nabla^2_a\varphi(w)
$$
$\mu$-a.s. in the strong topology of Hilbert-Schmidt operators.
\end{remarkk}

\noindent
Let us define, with $w=w_n^\bot+w_n=\pi^\bot_n(w)+\pi_n(w)$, the functional
$$
\delta_a(\pi_n\nabla\varphi)(w)=\sum_{i=1}^n(\nabla\varphi,e_i)_H\delta
e_i-\trace D_n^2(w_n^\bot+\cdot)_a(w_n)
$$
where the lower index represents the absolutely continuous component
of the measure on the Euclidean space $H_n$ spanned by the orthonormal
set $\{e_1,\ldots, e_n\}$. It follows then from Lemmas
\ref{section-lemma}, \ref{Bayes} and \ref{fin-compo} the following
relation:
\begin{eqnarray}
\label{ma_eqn}
\lefteqn{-\delta_a(\pi_n\nabla\varphi)(w_n^\bot,w_n)+\log\dett(I_{\R^n}+\pi_n\nabla^2\varphi\pi_n(w_n^\bot+\cdot)_a(w_n))-\half|\pi_n\nabla\varphi|^2}\\
&&=-\log
l_n(w_n^\bot+\cdot)\circ(w_n+\pi_n\nabla\varphi(w^\bot_n+w_n))\,,\nonumber
\end{eqnarray}
where 
$$
l_n(w_n^\bot+w_n)=\frac{L(w_n^\bot+w_n)}{E[L|\pi_n^\bot=w_n^\bot]}\,.
$$

\noindent
The following result gives the key information for the sequel:
\begin{theorem}
\label{RN-n}
Let $T_n(w)=w+\pi_n\nabla\varphi(w)$, then $T_n\mu$ is absolutely
continuous w.r.t. $\mu$ with the corresponding Radon-Nikodym
derivative
$$
M_n=\frac{dT_n\mu}{d\mu}=\frac{L}{E[L|\pi_n^\bot]}
$$
$\mu$-almost surely. Moreover $(M_n,n\geq 1)$ is uniformly integrable
(w.r.t. $\mu$).
\end{theorem}
\nproof
Let $g\in C_b(W)$, it follows from Lemma \ref{bayes_1} and Lemma
\ref{fin-compo} that
\beaa
E[g\circ T_n]&=&E[E[g\circ T_n|\pi_n^\bot]]\\
&=&\int E[g\circ T_n|w_n^\bot]\mu_n^\bot(dw_n^\bot)\\
&=&\int E\left[g\,\frac{L}{E[L|w_n^\bot]}\big|w_n^\bot\right]\mu_n^\bot(dw_n^\bot)\\
&=&\int E[g\,L|w_n^\bot]
 \frac{1}{E[L|w_n^\bot]} \mu_n^\bot(dw_n^\bot)\\
&=&E\left[g\frac{L}{E[L|\pi_n^\bot]}\right]\,.
\eeaa
Let us remark that, from the elementary properties of the conditional
expectation, $E[L|\pi_n^\bot]\neq 0$ on the set $\{L\neq 0\}$
$\mu$-almost surely. Finally, it follows from the martingale
convergence theorem $(M_n=\frac{L}{E[L|\pi_n^\bot]},n\geq 1)$ converges to $L$ almost
surely and $E[M_n]=1=E[L]$ for any $n\geq 1$, hence the uniform
integrability follows.
\nqed

The following theorem is the Monge-Amp\`ere equation satisfied by
$(\varphi,\psi)$:
\begin{theorem}
\label{MA_thm}
The sequence of Wiener functionals $(\delta_a\pi_n\nabla\varphi,n\geq 1)$
converges $\mu$-almost surely to a Wiener functional
$\calL^a\varphi$. Moreover the Wiener Jacobian defined by
$$
\La(\varphi)=\dett(I_H+\nabla^2_a\varphi)\exp\left[-\calL^a\varphi-\half|\nabla\varphi|_H^2\right]\,,
$$
satisfies the Monge-Amp\`ere equation:
\begin{equation}
\label{MA1_eqn}
L\circ T\,\La(\varphi)=1
\end{equation}
$\mu$-a.s., where $T=I_W+\nabla\varphi$. In particular, if
$\mu\{L>0\}=1$, then, for any $g\in C_b(W)$, we have the
Jacobi-Girsanov relation:
$$
\int_W g\circ T \La(\varphi)d\mu=\int_W gd\mu\,.
$$
\end{theorem}
\nproof
We shall use the notations of this section without further
recall. First note that, due to the uniform integrability of
$(M_n,n\geq 1)$ proven in Theorem \ref{RN-n}, an application of  Lusin
theorem implies that 
$$
\lim_{n\to\infty}\frac{L\circ T_n}{E[L|\pi_n^\bot]}=L\circ T
$$
$\mu$-a.s. It follows from Remark \ref{cv_rem} that 
$$
\lim_{n\to\infty}\dett(I_{\R^n}+\pi_n\nabla^2\varphi\pi_n(\pi_n^\bot(w)+\cdot)_a(\pi_n(w))=\dett(I+\nabla^2_a\varphi(w))
$$
$\mu$-a.s. As $\lim_n|\pi_n\nabla\varphi|_H^2=|\nabla\varphi|_H^2$
$\mu$-a.s., also, it follows from the equality (\ref{ma_eqn}), the sequence
$(\delta_a\pi_n\nabla\varphi,n\geq 1)$ converges $\mu$-almost surely
to a Wiener functional that we denote as $\calL^a\varphi$ and the
relation (\ref{MA1_eqn}) follows. To show the Jacobi-Girsanov relation,
for any $g\in C_b(W)$, we have, from (\ref{MA1_eqn}),
\beaa
\int g\circ T\La(\varphi)d\mu&=&\int g\circ T\frac{1}{L\circ T}d\mu\\
&=&\lim_{\eps\to 0}\int g\circ T\frac{1}{L\circ T+\eps}d\mu\\
&=&\lim_{\eps\to 0}\int g \frac{L}{L+\eps}d\mu\\
&=&\int_{\{L\neq 0\}}gd\mu=\int g d\mu\,,
\eeaa
since $\mu(\{L\neq 0\})=1$.
\nqed

\section{\bf{Regularity of the forward potential $\varphi$}}

\noindent
Let us show now the regularity of the forward Monge potential
$\varphi$: assume first 
that, we have reduced the problem to the case where 
everything is smooth using the approximation results that we have proven before. Let $\nu$ be the measure defined by $d\nu=e^{-f}d\mu$. The 
following  relation holds then true:
\begin{eqnarray}
-\log\nu(e^f)&=&\inf_\alpha\left(\int-f d\alpha+H(\alpha|\nu)\right) \label{1}\\
&=&\inf_U\left(\int-f\circ U d\nu+H(U\nu|\nu)\right)\,\label{2}.
\end{eqnarray}
It is important to remark that in the equation \ref{1}, the infimum is
taken over the set of probability measures and in the equation
\ref{2}, the infimum is taken over the perturbations of identity of
the form $U=I_W+u$ when $u$ runs in the set of the gradients of
$1$-convex functions, cf. \cite{fandu2}.
Moreover, denoting $\frac{dU\mu}{d\nu}$ by $l_U$, we have 
$$
(l_U\, e^{-f})\circ U\La_u=e^{-f}, 
$$
where $\La_u$ is the Gaussian Jacobian associated to
$U=I_W+u$. Therefore $\log l_U\circ U=f\circ U-\log \La_u-f$ and we
get
\beaa
H(U\nu|\nu)&=&\int (f\circ U-\log \La_u-f)d\nu\\
&=&\int (f\circ U-\log \La_u-f)e^{-f}d\mu\,.
\eeaa
Consequently
\beaa
-\log\nu(e^f)&=&\inf_U\left(\int-f\circ U d\nu+\int f\circ U
  e^{-f}d\mu-\int(f+\log\La_u)e^{-f}d\mu\right)\\
&=&\inf_U\left(-\int fe^{-f}d\mu-\int\log\La_ue^{-f}d\mu\right)\\
&=&\inf_U J_b(U)\,.
\eeaa
We know that the above  infimum is attained at $S=T^{-1}=I_W+\nabla\psi$, hence we
should have 
$$
J_b'(S)\cdot\xi=\frac{d}{d\la}J_b(S+\la \xi)|_{\la=0}=0
$$
for any smooth $\xi:W\to H$ such that $\|\nabla\xi\|_2\in
L^\infty(\mu)$. A similar calculation as performed before implies that
\beaa
\frac{d}{d\la}J_b(S+\la \xi)|_{\la=0}&=&\frac{d}{d\la}\left(\int
-\log\La_{S+\la\xi} d\nu\right)\Big|_{\la=0}\\
&=&\int
\left[-\trace(((I+\nabla^2\psi)^{-1}-I)\cdot\nabla\xi)+\delta\xi+(\nabla\psi,\xi)\right]d\nu\\
&=&0
\eeaa
for any $\xi$ as above. Consequently we have
\begin{theorem}
\label{d-potential}
The backward Monge potential satisfies the relation
$$
\delta_\nu((I+\nabla^2\psi)^{-1}-I)=\nabla\psi+\nabla f\,,
$$
where $\delta_\nu$ denotes the adjoint of $\nabla$ w.r.t. the measure
$\nu$.
\end{theorem}

\noindent
We need a couple of techical results:
\begin{lemma}
\label{lemma_1}
Let $\xi:W\to H$ be a smooth vector field, then the following results
hold true:
\begin{enumerate}
\item $\delta_\nu\xi=\delta\xi+(\nabla f,\xi)_H$.
\item For any $h\in H$, 
$$
E_\nu[(\delta_\nu h)^2]=E_\nu[|h|^2+(\nabla^2f,h\otimes h)]\,.
$$
\item For any $h\in H$ and smooth $\alpha:W\to \R$, 
$$
E_\nu[\alpha(\delta_\nu h)^2]=E_\nu\left[(\alpha
I_H+\nabla^2\alpha+\alpha\nabla^2f,h\otimes h)_{H^{\otimes
    2}}\right]\,.
$$
\end{enumerate}
\end{lemma}

\noindent
\begin{lemma}
\label{lemma_2}
For any smooth $\xi:W\to H$, we have
$$
E_\nu[(\delta_\nu\xi)^2]=E_\nu\left[(I_H+\nabla^2f,\xi\otimes
  \xi)_{H^{\otimes 2}}+\trace(\nabla\xi\cdot\nabla\xi)\right]\,.
$$
\end{lemma}
\nproof
By the definition of $\delta_\nu$, we have
$$
E_\nu[\alpha(\delta_\nu h)^2]=E_\nu[|\xi|^2+(\xi,\delta\otimes
\nabla\xi)+(\xi,\nabla(\nabla f,\xi))]\,.
$$
Besides
$(\xi,\delta\otimes\nabla\xi)=\delta\nabla_\xi\xi+\trace(\nabla\xi\cdot\nabla\xi)$
(cf.\cite{BOOK}). Hence
$$
E_\nu[\alpha(\delta_\nu h)^2]=E_\nu[|\xi|^2+\delta\nabla_\xi\xi+\trace(\nabla\xi\cdot\nabla\xi)
+\delta_\nu\xi\,(\nabla f,\xi))]\,.
$$
We also have
$$
\delta\nabla_\xi\xi=\delta_\nu\nabla_\xi\xi-(\nabla f,\nabla_\xi,\xi)\,.
$$
Substituting this expression in the above calculation  gives
\beaa
E_\nu[(\delta_\nu\xi)^2]&=&E_\nu[|\xi|_H^2+\delta_nu(\nabla_\xi\xi)-(\nabla
  f,\nabla_\xi\xi)_H+\trace(\nabla\xi\cdot\nabla\xi)\\
&&+\delta_\nu\xi\,(\nabla f,\xi)_H]\\
&=&E_\nu\left[|\xi|_H^2-(\nabla
  f,\nabla_\xi\xi)_H+\trace(\nabla\xi\cdot\nabla\xi)+(\xi,\nabla(\nabla
  f,\xi)_H)_H\right]\\
&=&E_\nu\left[|\xi|_H^2+\trace(\nabla\xi\cdot\nabla\xi)-(\nabla
  f,\nabla_\xi\xi)_H+(\nabla_\xi\nabla f,\xi)_H+(\nabla
  f,\nabla_\xi\xi)_H\right]\\
&=&E_\nu\left[|\xi|_H^2+\trace(\nabla\xi\cdot\nabla\xi)+(\nabla_\xi\nabla
  f,\xi)_H\right]\\
&=&E_\nu\left[|\xi|_H^2+\trace(\nabla\xi\cdot\nabla\xi)+(\nabla^2f,\xi\otimes\xi)_{H^{\otimes
      2}}\right]\,.
\eeaa
\nqed

\noindent
The following theorem extends a result of Caffarelli \cite{Caf1}, from
log-concave densities to $1-\eps$-log-concave densities with a
different proof:
\begin{theorem}
\label{forward_thm}
Assume that $f\in L^{p}(\mu)$ for some $p>1$, satisfying $E[|\nabla
f|_H^2e^{-f}]<\infty$. Assume moreover that it is $(1-\eps)$-convex for
some $\eps>0$, in the sense that the mapping
$$
h\to \frac{1-\eps}{2}|h|_H^2+f(w+h)
$$
is a convex map from the Cameron-Martin space $H$  to $L^0(\mu)$ (i.e., the
equivalence class of real-valued Wiener functionals under the topology
of convergence in probability). Then the forward Monge potential $\varphi$ belongs to
the Gaussian Sobolev space $\DD_{2,2}$.
\end{theorem}
\nproof
let $f_n, n\geq 1$ be defined as
$e^{-f_n}=P_{1/n}E[e^{-f}|V_n]$. Since $f_n$ is a smooth, $1$-convex
function, the corresponding forward potential $\varphi_n$ is also
smooth from the classical finite dimensional results (cf. \cite{Caf},
\cite{Vil}). Let $d\nu_n=e^{-f_n}d\mu$, then we have, from Theorem
\ref{d-potential}
$$
\delta_{\nu_n}((I_H+\nabla^2\psi_n)^{-1}-I_H)=\nabla\psi_n+\nabla
f_n\,.
$$
From Lemma \ref{lemma_2} and denoting $(I_H+\nabla^2\psi_n)^{-1}$ by
$M_n$, we get
\beaa
E_{\nu_n}\left[|\delta_{\nu_n}((I_H+\nabla^2\psi_n)^{-1}-I_H)|_H^2\right]&=&\sum_{k=1}^\infty 
E_{\nu_n}[(\delta_{\nu_n}(M_n-I)(e_k))^2]\\
&=&\sum_{k=1}^\infty E_{\nu_n}[(I_H+\nabla^2f_n, (M_n-I_H)e_k\otimes
(M_n-I_H)e_k)]\\
&&\hspace{1cm}+\sum_{k=1}^\infty E_{\nu_n}[\trace(\nabla (M_ne_k)\cdot(\nabla
M_ne_k))]\,.
\eeaa
Since, due to the $(1-\eps)$-convexity of $f$ and as the second terms at the right of the second line is positive, we obtain
$$
E_{\nu_n}[|\delta_{\nu_n}(M_n-I_H)|_H^2]\geq \eps\sum_{k=1}^\infty
E_{\nu_n}[|(M_n-I_H)e_k|_H^2]\,.
$$
Hence
$$
\eps E_{\nu_n}\left[\|(I_H+\nabla^2\psi_n)^{-1}-I_H\|_2^2\right]\leq 2
E_{\nu_n}\left[|\nabla\psi_n|_H^2+|\nabla f_n|_H^2\right]\,,
$$
but $E_{\nu_n}[|\nabla\psi_n|_H^2]=E[|\nabla\varphi_n|_H^2]$ and 
\beaa
E_{\nu_n}\left[\|(I_H+\nabla^2\psi_n)^{-1}-I_H\|_2^2\right]&=&
E\left[\|(I_H+\nabla^2\psi_n)^{-1}\circ T_n-I_H\|_2^2\right]\\
&=&E[\|\nabla^2\varphi_n\|_2^2]\,.
\eeaa
We also have 
\beaa
E_{\nu_n}\left[|\nabla f_n|_H^2\right]&=&4 E[|\nabla e^{-f_n/2}|^2]\\
&=&4E\left[\frac{1}{e^{-f_n}}\,|\nabla f_ne^{-f_n}|^2\right]\\
&=&4 E\left[\frac{1}{e^{-f_n}} |\nabla e^{-f_n}|^2\right]\\
&=& 4 E\left[\frac{1}{e^{-f_n}} |\nabla P_{1/n}E[ e^{-f}|V_n]|^2\right]\\
&\leq&4e^{-1/n}E\left[\frac{1}{e^{-f_n}} | P_{1/n}E[ \nabla
e^{-f}|V_n]|^2\right]\\
&\leq&4e^{-1/n}E\left[\frac{1}{e^{-f_n}} P_{1/n}(E[|\nabla
f|^2e^{-f}|V_n])P_{1/n}E[e^{-f}|V_n]\right]\\
&=&4e^{-1/n} E[P_{1/n}E[|\nabla f|^2e^{-f}|V_n]]\\
&=&4e^{-1/n} E[|\nabla f|^2e^{-f}]\,.
\eeaa
Consequently we get
$$
\eps E[\|\nabla^2\varphi_n\|_2^2]\leq 2E[|\nabla\varphi_n|_H^2]+8E[|\nabla
f|^2e^{-f}] 
$$
and the claim follows by taking the limit at the r.h.s. and he limit
inferior at the l.h.s. even with an explicit bound:
$$
\eps E[\|\nabla^2\varphi\|_2^2]\leq 2E[|\nabla\varphi|_H^2]+8E[|\nabla
f|^2e^{-f}] \,.
$$
\nqed

\noindent
The next corollary  which  follows from Theorem \ref{psi_cvg_thm} and from
Lemma \ref{lemma_2}, is about the regularity of the dual potential $\psi$:
\begin{corollary}
Assume that $E[\|\nabla^2f\|^2_\infty e^{-f}]<\infty$, where
$\|\cdot\|_\infty$ denotes the operator norm on $H$. Then
$(\delta_\nu\circ \nabla) \psi=\calL_\nu\psi$ belongs to $L^2(\nu)$. 
\end{corollary}

\vspace{2cm}
\footnotesize{
\noindent
A. S. \"Ust\"unel, Bilkent University, Math. Dept., Ankara, Turkey\\
ustunel@fen.bilkent.edu.tr}

\end{document}